\documentclass[11pt]{article}
\usepackage{amssymb,latexsym}
\usepackage{epsfig}
\usepackage{eufrak}
\usepackage{amsmath}
\usepackage{mathrsfs}
\usepackage{color}

\newtheorem{theorem}{Theorem}
\newtheorem{lemma}{Lemma}

\newcommand{\be}{\begin{equation}}
\newcommand{\ee}{\end{equation}}
\newcommand{\bee}{\begin{eqnarray*}}
\newcommand{\eee}{\end{eqnarray*}}
\newcommand{\bel}{\begin{eqnarray}}
\newcommand{\eel}{\end{eqnarray}}
\newcommand{\bec}{\begin{cases}}
\newcommand{\eec}{\end{cases}}
\newcommand{\bem}{\begin{bmatrix}}
\newcommand{\eem}{\end{bmatrix}}

\newcommand{\la}{\label}
\newcommand{\li}{\left}
\newcommand{\ri}{\right}

\newcommand{\ovl}{\overline}
\newcommand{\udl}{\underline}

\newcommand{\lc}{\lceil}
\newcommand{\rc}{\rceil}
\newcommand{\lf}{\lfloor}
\newcommand{\rf}{\rfloor}

\newcommand{\ep}{\epsilon}
\newcommand{\vep}{\varepsilon}

\newcommand{\de}{\delta}

\newcommand{\ga}{\gamma}

\newcommand{\se}{\theta}
\newcommand{\Se}{\Theta}

\newcommand{\ze}{\zeta}
\newcommand{\al}{\alpha}

\newcommand{\ro}{\rho}

\newcommand{\om}{\omega}
\newcommand{\Om}{\Omega}

\newcommand{\f}{\frac}
\newcommand{\sq}{\sqrt}
\newcommand{\cd}{\cdots}

\newcommand{\qu}{\quad}
\newcommand{\qqu}{\qquad}
\newcommand{\fa}{\forall}

\newcommand{\mscr}{\mathscr}
\newcommand{\mcal}{\mathcal}
\newcommand{\mbf}{\mathbf}
\newcommand{\bb}{\mathbb}

\newcommand{\wh}{\widehat}
\newcommand{\wt}{\widetilde}
\newcommand{\mrm}{\mathrm}
\newcommand{\bs}{\boldsymbol}

\newcommand{\sh}{\slash}

\newcommand{\tx}{\text}

\newcommand{\iy}{\infty}

\newcommand{\leu}{\subseteq}

\newcommand{\pa}{\partial}

\newcommand{\bed}{\begin{description}}
\newcommand{\eed}{\end{description}}
\newcommand{\bei}{\begin{itemize}}
\newcommand{\eei}{\end{itemize}}
\newcommand{\ben}{\begin{enumerate}}
\newcommand{\een}{\end{enumerate}}
\newcommand{\bib}{\bibitem}
\newcommand{\beL}{\begin{lemma}}
\newcommand{\eeL}{\end{lemma}}
\newcommand{\beT}{\begin{theorem}}
\newcommand{\eeT}{\end{theorem}}
\newcommand{\sect}{\section}

\newcommand{\bpf}{\begin{pf}}
\newcommand{\epf}{\end{pf}}
\newcommand{\bsk}{\bigskip}

\newcommand{\pfbox}{\hfill\mbox{$\Box$}}

\newenvironment{pf}{\paragraph*{Proof{\rm.}}}{\pfbox\bigskip}

\begin{document}

\title{{\bf Multistage Estimation of Bounded-Variable Means}
\thanks{The author had been previously working with Louisiana State University at Baton Rouge, LA 70803, USA,
and is now with Department of Electrical Engineering, Southern
University and A\&M College, Baton Rouge, LA 70813, USA; Email:
chenxinjia@gmail.com}}

\author{Xinjia Chen}

\date{September 2008}

\maketitle

\begin{abstract}

In this paper, we develop a multistage approach for estimating the
mean of a bounded variable.  We first focus on the multistage
estimation of a binomial parameter and then generalize the
estimation methods to the case of general bounded random variables.
A fundamental connection between a binomial parameter and the mean
of a bounded variable is established.  Our multistage estimation
methods rigorously guarantee prescribed levels of precision and
confidence.

\end{abstract}

\sect{Introduction}

The estimation of the means of bounded random variables finds
numerous applications in various fields of sciences and engineering.
In particular, Bernoulli random variables constitute an extremely
important class of bounded variables, since the ubiquitous problem
of estimating the probability of an event can be formulated as the
estimation of the mean of a Bernoulli variable. In many
applications, one needs to estimate a quantity $\mu$ which can be
bounded in $[0,1]$ after proper operations of scaling and
translation.  A typical approach is to design an experiment that
produces a random variable $Z$ distributed in $[0, 1]$ with
expectation $\mu$, run the experiment independently a number of
times, and use the average of the outcomes as the estimate
\cite{KLM}. This technique, referred to as {\it Monte Carlo} method,
has been applied to tackle a wide range of difficult problems.

Since the estimator of the mean of $Z$ is obtained from finite
samples of $Z$ and is thus of random nature, for the estimator to be
useful, it is necessary to ensure with a sufficiently high
confidence that the estimation error is within certain margin.  The
well known  Chernoff-Hoeffding bound \cite{Chernoff}
\cite{Hoeffding} asserts that if the sample size is fixed and is
greater than $\f{\ln \f{2}{\de} } {2 \ep^2 }$, then, with
probability at least $1 - \de$, the sample mean approximates $\mu$
with absolute error $\ep$.   The problem with Chernoff-Hoeffding
bound is that the resultant sample size can be extremely
conservative if the value of $\mu$ is close to zero or one. In the
case that $\mu$ is small, it is more reasonable to seek an $(\vep,
\de)$ approximation for $\mu$ in the sense that the relative error
of the estimator is within a margin of relative error $\vep$ with
probability at least $1 - \de$.  Since
 the mean value $\mu$ is exactly what we want to estimate, it is usually not easy to obtain reasonably tight lower bound
 for $\mu$.  For a sampling scheme with fixed sample size, a loose lower bound of $\mu$ can lead to a very
 conservative sample size.  For the most difficult and important case that no positive lower bound of $\mu$ is available,
it is not possible to guarantee prescribed relative precision and
confidence level by a sampling scheme with a fixed sample size.
This forces us to look at sampling methods with random sample sizes.

The estimation techniques based on sampling schemes without fixed
sample sizes have formed a rich branch of modern statistics under
the heading of {\it sequential estimation}.  Wald provided a brief
introduction to this area in his seminal book \cite{Wald}.  Ghosh et
al. offered a comprehensive exposition in \cite{Gosh}. In
particular, Nadas proposed in \cite{Naddas} a sequential sampling
scheme for estimating mean values with relative precision.  Nadas's
sequential method requires no specific information on the mean value
to be estimated. However, his sampling scheme is of asymptotic
nature. The confidence requirement is guaranteed only as the margin
of relative error $\vep$ tends to $0$, which implies that the actual
sample size has to be infinity. This drawback severely circumvents
the application of his sampling scheme.

In this paper, we revisit  the sequential estimation of means of
random variables bounded in $[0,1]$.  To overcome the limitations of
existing methods, we have developed a new class of multistage
sampling schemes. Our sampling schemes require no information of the
unknown parameters and guarantees prescribed levels of precision and
confidence.  The remainder of the paper is organized as follows.
Section 2 is devoted to the multistage estimation of a binomial
parameter. In Section 3, we generalize the estimation methods of a
binomial parameter to the mean of a bounded variable.  In Section 4,
we establish a link between a binomial parameter and the mean of a
bounded variable. We demonstrate that the estimation methods for
estimating a binomial parameter can be easily applied to the
estimation of the mean of a bounded variable by virtue of this link.
Section 5 is the conclusion. All proofs are given in the Appendices.

Throughout this paper, we shall use the following notations. The
expectation of a random variable is denoted by $\bb{E}[.]$. The set
of integers is denoted by $\bb{Z}$. The set of positive integers is
denoted by $\bb{N}$. The ceiling function and floor function are
denoted respectively by $\lc . \rc$ and $\lf . \rf$ (i.e., $\lc x
\rc$ represents the smallest integer no less than $x$; $\lf x \rf$
represents the largest integer no greater than $x$).  The notation
$\mrm{sgn}(x)$ denotes the sign function which assumes value $1$ for
$x > 0$, value $0$ for $x = 0$, and value $-1$ for $x < 0$.  We use
the notation $\Pr \{ . \mid \se \}$ to indicate that the associated
random samples $X_1, X_2, \cd$ are parameterized by $\se$. The
parameter $\se$ in $\Pr \{ . \mid \se \}$  may be dropped whenever
this can be done without introducing confusion. The other notations
will be made clear as we proceed.

\sect{Estimation of Binomial Parameters}

Let $X$ be a Bernoulli random variable defined in a probability
space $(\Om, \mscr{F}, \Pr )$ such that $\Pr \{ X = 1 \} = 1 - \Pr
\{ X = 0 \} = p \in (0, 1)$.  It is a frequent problem to estimate
the binomial parameter $p$ based on a sequence of i.i.d. random
samples $X_1, X_2, \cd$ of $X$.  To solve this problem, we shall
develop multistage sampling schemes of the following general
structure. The sampling process is divided into $s$ stages. The
continuation or termination of sampling is determined by decision
variables.  For each stage with index $\ell$, a decision variable
$\bs{D}_\ell = \mscr{D}_\ell (X_1, \cd, X_{\mbf{n}_\ell})$ is
defined based on samples $X_1, \cd, X_{\mbf{n}_\ell}$, where
$\mbf{n}_\ell$ is the number of samples available at the $\ell$-th
stage.  It should be noted that $\mbf{n}_\ell$ can be a random
number, depending on specific sampling schemes. The decision
variable $\bs{D}_\ell$ assumes only two possible values $0, \; 1$
with the notion that the sampling is continued until $\bs{D}_\ell =
1$ for some $\ell \in \{1, \cd, s\}$. Since the sampling must be
terminated at or before the $s$-th stage, it is required that
$\bs{D}_s = 1$. For simplicity of notations, we also define
$\bs{D}_\ell = 0$ for $\ell = 0$ throughout the remainder of the
paper.

\subsection{Control of Absolute Error}

In many situations, it is desirable to construct an estimator for
$p$ with guaranteed absolute precision and confidence level.  For
this purpose, we have

\beT \la{coverage_abs} Let $0 < \vep < \f{1}{2}, \; 0 < \de < 1, \;
\ze > 0$ and $\ro > 0$. Let $n_1 < n_2 < \cd < n_s$ be the ascending
arrangement of all distinct elements of {\small $\li \{ \li \lc \li
( \f{24 \vep - 16 \vep^2}{9} \ri )^{ 1 - \f{i}{\tau} } \f{\ln
\f{1}{\ze \de}}{2 \vep^2} \ri \rc : i = 0, 1, \cd, \tau \ri \}$}
with {\small $\tau = \li \lc \f{ \ln \f{9}{24 \vep - 16 \vep^2} } {
\ln (1 + \ro) } \ri \rc$.}    For $\ell = 1, \cd, s$, define $K_\ell
= \sum_{i = 1}^{n_\ell} X_i, \; \wh{\bs{p}}_\ell = \f{ K_\ell
}{n_\ell}$  and $\bs{D}_\ell$ such that $\bs{D}_\ell = 1$ if {\small
$\li ( \li | \wh{\bs{p}}_\ell - \f{1}{2} \ri | - \f{2 \vep }{3} \ri
)^2 \geq \f{1}{4} + \f{ \vep^2 n_\ell } {2 \ln (\ze \de) }$}; and
$\bs{D}_\ell = 0$ otherwise.  Suppose the stopping rule is that
sampling is continued until $\bs{D}_\ell = 1$ for some $\ell \in
\{1, \cd, s\}$. Define $\bs{\wh{p}} = \f{\sum_{i=1}^{\mathbf{n}}
X_i}{\mathbf{n}}$ where $\mathbf{n}$ is the sample size when the
sampling is terminated. Define {\small  \[ \mscr{Q}^+ =
\bigcup_{\ell = 1}^s \li \{ \f{k}{n_\ell} + \vep \in \li (0,
\f{1}{2} \ri ) : k \in \bb{Z} \ri \} \bigcup \li \{ \f{1}{2} \ri \},
\qqu \mscr{Q}^- = \bigcup_{\ell = 1}^s \li \{ \f{k}{n_\ell} - \vep
\in \li (0, \f{1}{2} \ri ) : k \in \bb{Z} \ri \} \bigcup \li \{
\f{1}{2} \ri \}.
\]}
Then, a sufficient condition to guarantee $\Pr \li \{ \li |
\bs{\wh{p}} - p \ri | < \vep \mid p \ri \}
> 1 - \de$ for any $p \in (0, 1)$ is that {\small  \bel &  &
\sum_{\ell = 1}^s \Pr \{ \wh{\bs{p}}_\ell \geq p + \vep, \;
\bs{D}_{\ell - 1} = 0, \; \bs{D}_\ell = 1 \mid p \} < \f{\de}{2}
\qqu \fa p \in
\mscr{Q}^-, \la{2D1}\\
&  & \sum_{\ell = 1}^s \Pr \{ \wh{\bs{p}}_\ell  \leq p - \vep, \;
\bs{D}_{\ell - 1} = 0, \; \bs{D}_\ell = 1 \mid p \} < \f{\de}{2}
\qqu \fa p \in \mscr{Q}^+ \la{2D2}
 \eel}
where both (\ref{2D1}) and (\ref{2D2}) are satisfied if $0 < \ze <
\f{1}{2(\tau + 1)}$.

\eeT

\subsection{Control of Absolute and Relative Errors}

To construct an estimator satisfying a mixed criterion in terms of
absolute and relative errors with a prescribed confidence level, we
have

\beT \la{coverage_mixed}  Let $0 < \de < 1, \; \ze > 0$ and $\ro >
0$. Let $\vep_a$ and $\vep_r$ be positive numbers such that $0 <
\vep_a < \f{3}{8}$ and {\small $\f{6 \vep_a}{3 - 2 \vep_a } < \vep_r
< 1$}.  Let $n_1 < n_2 < \cd < n_s$ be the ascending arrangement of
all distinct elements of {\small $\li \{ \li \lc \li [ \f{3}{2} \li
( \f{1}{\vep_a} - \f{1}{\vep_r} - \f{1}{3} \ri ) \ri ]^{ \f{i } {
\tau } } \f{ 4(3 + \vep_r) } {9 \vep_r } \ln \f{1} { \ze \de } \ri
\rc: i = 0, 1, \cd, \tau \ri \}$ } with {\small $\tau = \li \lc \f{
\ln \li [ \f{3}{2} \li ( \f{1}{\vep_a} - \f{1}{\vep_r} - \f{1}{3}
\ri ) \ri ] } { \ln (1 + \ro) } \ri \rc$. }   For $\ell = 1, \cd,
s$, define {\small $K_\ell = \sum_{i = 1}^{n_\ell} X_i, \;
\wh{\bs{p}}_\ell = \f{K_\ell}{n_\ell},$} {\small
\[ \bs{D}_\ell = \bec 0 & \mrm{for} \; \f{1}{2} - \f{2}{3} \vep_a -
\sq{ \f{1}{4} + \f{ n_\ell \vep_a^2 } {2 \ln (\ze \de) } } <
\wh{\bs{p}}_\ell < \f{ 6(1 - \vep_r) (3 - \vep_r) \ln (\ze
\de) } { 2 (3 - \vep_r)^2 \ln (\ze \de) - 9 n_\ell \vep_r^2} \; \mrm{or}\\
  &  \qu \;\; \f{1}{2} + \f{2}{3} \vep_a - \sq{ \f{1}{4} + \f{
n_\ell \vep_a^2 } {2 \ln (\ze \de) } } < \wh{\bs{p}}_\ell < \f{ 6(1
+ \vep_r) (3 + \vep_r) \ln (\ze \de) } { 2 (3 +
\vep_r)^2 \ln (\ze \de) - 9 n_\ell  \vep_r^2},\\
1 & \mrm{else}
 \eec
\]}
for $\ell = 1, \cd, s - 1$ and $\bs{D}_s = 1$.  Suppose the stopping
rule is that sampling is continued until {\small $\bs{D}_\ell = 1$}
for some $\ell \in \{1, \cd, s\}$. Let {\small $\wh{\bs{p}} =
\f{\sum_{i=1}^{\mathbf{n}} X_i} {\mathbf{n}}$} where $\mathbf{n}$ is
the sample size when the sampling is terminated. Define {\small
$p^\star = \f{\vep_a}{\vep_r}$} and {\small
\[ \mscr{Q}_a^+ = \bigcup_{\ell = 1}^s \li \{ \f{k}{n_\ell} + \vep_a
\in \li (0, p^\star \ri ) : k \in \bb{Z} \ri \} \cup \li \{ p^\star
\ri \}, \qqu \mscr{Q}_a^- = \bigcup_{\ell = 1}^s \li \{
\f{k}{n_\ell} - \vep_a \in \li (0, p^\star \ri ) : k \in \bb{Z} \ri
\} \cup \li \{ p^\star \ri \},
\]
\[
\qqu \mscr{Q}_r^+ = \bigcup_{\ell = 1}^s \li \{ \f{k}{n_\ell (1 +
\vep_r)} \in \li (p^\star, 1 \ri ) : k \in \bb{Z} \ri \}, \qu \qqu
\qu \mscr{Q}_r^- = \bigcup_{\ell = 1}^s \li \{ \f{k}{n_\ell (1 -
\vep_r)} \in \li (p^\star, 1 \ri ) : k \in \bb{Z} \ri \}.\qqu \qqu
\qqu\qqu
\]}
Then, {\small $\Pr \li \{ \li | \wh{\bs{p}} - p
 \ri | < \vep_a  \; \mrm{or} \; \li |
\f{\wh{\bs{p}} - p } {p }
 \ri | < \vep_r \mid p \ri \} > 1 - \de$
 }
for any $p \in (0, 1)$  provided that {\small \bel &   & \sum_{\ell
= 1}^s \Pr \{ \wh{\bs{p}}_\ell \geq p + \vep_a, \; \bs{D}_{\ell - 1}
= 0, \; \bs{D}_\ell = 1 \mid p \} < \f{\de}{2} \qqu \fa p \in
\mscr{Q}_a^-, \la{mix1}\\
&   & \sum_{\ell = 1}^s \Pr \{ \wh{\bs{p}}_\ell \leq p - \vep_a, \;
\bs{D}_{\ell - 1} = 0, \; \bs{D}_\ell = 1 \mid p \} < \f{\de}{2}
\qqu \fa p \in \mscr{Q}_a^+, \la{mix2}\\
&   & \sum_{\ell = 1}^s \Pr \{ \wh{\bs{p}}_\ell \geq p (1 + \vep_r),
\; \bs{D}_{\ell - 1} = 0, \; \bs{D}_\ell = 1 \mid p \} < \f{\de}{2}
\qqu \fa p \in \mscr{Q}_r^+, \la{mix3}\\
&   &  \sum_{\ell = 1}^s \Pr \{ \wh{\bs{p}}_\ell \leq p (1 -
\vep_r), \; \bs{D}_{\ell - 1} = 0, \; \bs{D}_\ell = 1 \mid p \} <
\f{\de}{2} \qqu \fa p \in \mscr{Q}_r^- \la{mix4}
 \eel}  where these conditions are satisfied for $0 < \ze < \f{1}{2 (\tau + 1) }$ .  \eeT

\subsection{Control of Relative Error}

In many situations, it is desirable to design a sampling scheme to
estimate $p$ such that the estimator satisfies a relative error
criterion with a prescribed confidence level.  By virtue of the
function {\small \[ g(\vep, \ga) = 1 - \sum_{i= 0}^{\ga - 1}
\f{1}{i!} \li ( \f{\ga}{ 1 + \vep} \ri )^i \exp \li ( - \f{\ga}{ 1 +
\vep} \ri ) + \sum_{i= 0}^{\ga - 1} \f{1}{i!}  \li ( \f{\ga}{ 1 -
\vep} \ri )^i \exp \li ( - \f{\ga}{ 1 - \vep} \ri ),
\]}
we have developed a simple sampling scheme as described by the
following theorem.

\beT \la{coverage_rev_a} Let $0 < \vep < 1, \; 0 < \de < 1, \; \ze >
0$ and $\ro > 0$. Let $\ga_1 < \ga_2 < \cd < \ga_s$ be the ascending
arrangement of all distinct elements of {\small $\li \{ \li \lc \li
[  \f{3}{2} \li ( \f{1}{\vep} + 1 \ri ) \ri ]^{ \f{i} { \tau } } \f{
4(3 + \vep) } {9 \vep } \ln \f{1} { \ze \de } \ri \rc: i = 0, 1,
\cd, \tau \ri \}$ } with {\small $\tau = \li \lc \f{ \ln \li [
\f{3}{2} \li ( \f{1}{\vep} + 1 \ri )  \ri ] } { \ln (1 + \ro) } \ri
\rc$.}   Let $\wh{\bs{p}}_\ell = \f{ \sum_{i = 1}^{\mathbf{n}_\ell}
X_i } { \mathbf{n}_\ell }$ where $\mathbf{n}_\ell$ is the minimum
number of samples such that $\sum_{i = 1}^{\mathbf{n}_\ell} X_i =
\ga_\ell$. For $\ell = 1, \cd, s$,  define $\bs{D}_\ell$ such that
$\bs{D}_\ell = 1$ if {\small $\wh{\bs{p}}_\ell \geq 1 + \f{ 2
\vep}{3 + \vep} + \f{ 9 \vep^2 \ga_\ell } { 2 (3 + \vep)^2 \ln (\ze
\de)}$}; and $\bs{D}_\ell = 0$ otherwise. Suppose the stopping rule
is that sampling is continued until $\bs{D}_\ell = 1$ for some $\ell
\in \{1, \cd, s\}$. Define estimator $\wh{\bs{p}} =
\f{\sum_{i=1}^{\mathbf{n}} X_i}{\mathbf{n}}$ where $\mathbf{n}$ is
the sample size when the sampling is terminated. Then, {\small $\Pr
\li \{ \li | \f{ \wh{\bs{p}} - p } { p }  \ri | \leq \vep \mid p \ri
\} \geq 1 - \de$} for any $p \in (0, 1)$ provided that $\ze > 0$ is
sufficiently small to guarantee $g(\vep, \ga_s) < \de$ and
 {\small \bel &  & \ln ( \ze \de )
< \li [ \f{ \li ( 1 + \vep  +  \sq{ 1 + 4 \vep + \vep^2 }  \ri )^2}
{ 4 \vep^2 } + \f{1}{2} \ri ] \li [ \f{\vep}{ 1 + \vep}  - \ln (1 +
\vep) \ri ], \la{cona}\\
&  & \sum_{\ell = 1}^s \Pr \{ \wh{\bs{p}}_\ell \leq (1 - \vep) p, \;
\bs{D}_{\ell - 1} = 0, \; \bs{D}_\ell = 1 \mid p \} \leq \f{\de}{2}
\qqu \fa p \in \mscr{Q}_r^-, \la{rev1}\\
&  & \sum_{\ell = 1}^s \Pr \{ \wh{\bs{p}}_\ell \geq (1 + \vep) p, \;
\bs{D}_{\ell - 1} = 0, \; \bs{D}_\ell = 1 \mid p \} \leq \f{\de}{2}
\qqu \fa p \in \mscr{Q}_r^+  \la{rev2} \eel} where {\small
$\mscr{Q}_r^+ = \bigcup_{\ell = 1}^s \li \{\f{\ga_\ell}{m (1 + \vep)
} \in (p^*, 1) : m \in \bb{N} \ri \} $} and {\small $\mscr{Q}_r^- =
\bigcup_{\ell = 1}^s \li \{\f{\ga_\ell}{m (1 - \vep) } \in (p^*, 1)
: m \in \bb{N} \ri \}$} with $p^* \in (0, z_{s-1})$ denoting the
unique number satisfying \[ g(\vep, \ga_s)  + \sum_{\ell = 1}^{s -
1} \exp \li ( \f{\ga_\ell}{z_\ell} \f{ (p^* - z_\ell)^2 } {2 \li (
\f{2 p^*}{3} + \f{z_\ell}{3} \ri ) \li ( \f{2 p^*}{3} +
\f{z_\ell}{3} - 1 \ri ) } \ri ) = \de \]
 where {\small $z_\ell =  1 + \f{ 2 \vep}{3 + \vep} + \f{ 9 \vep^2 \ga_\ell } { 2 (3
+ \vep)^2 \ln (\ze \de)}$} for $\ell = 1, \cd, s - 1$.
 \eeT

 In this section, we have proposed a multistage inverse sampling plan
for estimating a binomial parameter, $p$, with relative precision.
In some situations, the cost of sampling operation may be high since
samples are obtained one by one when inverse sampling is involved.
In view of this fact, it is desirable to develop multistage
estimation methods without using inverse sampling.  For this
purpose, we have

\beT \la{noinverse} Let $0 < \vep < 1, \; 0 < \de < 1$ and $\ze
> 0 $.  Let $\tau$ be a positive integer. For $\ell = 1, 2, \cd$, let $\wh{\bs{p}}_\ell = \f{ \sum_{i =
1}^{n_\ell} X_i } { n_\ell }$, where $n_\ell$ is deterministic and
stands for the sample size at the $\ell$-th stage. For $\ell = 1, 2,
\cd$, define $\bs{D}_\ell$ such that $\bs{D}_\ell = 1$ if {\small
$\wh{\bs{p}}_\ell \geq \f{ 6(1 + \vep) (3 + \vep) \ln (\ze \de_\ell)
} { 2 (3 + \vep)^2 \ln (\ze \de_\ell) - 9 n_\ell \vep^2}$}; and
$\bs{D}_\ell = 0$ otherwise, where $\de_\ell = \de$ for $1 \leq \ell
\leq \tau$ and $\de_\ell = \de 2^{\tau - \ell}$ for $\ell > \tau$.
Suppose the stopping rule is that sampling is continued until
$\bs{D}_\ell = 1$ for some stage with index $\ell$. Define estimator
$\wh{\bs{p}} = \wh{\bs{p}}_{\bs{l}}$, where $\bs{l}$ is the index of
stage at which the sampling is terminated. Then, $\Pr \{ \bs{l} <
\iy \} = 1$ and {\small $\Pr \li \{ \li | \f{ \wh{\bs{p}} - p } { p
} \ri | \leq \vep \mid p \ri \} \geq 1 - \de$} for any $p \in (0,
1)$ provided that $2 (\tau + 1 ) \ze \leq 1$ and $\inf_{\ell > 0}
\f{n_{\ell + 1}}{n_\ell} > 0$.

\eeT

\subsection{Fixed-width Confidence Intervals}

In some literature, the estimation of $p$ has been formulated as the
problem of constructing a fixed-width confidence interval $(\bs{L},
\bs{U})$ such that $\bs{U} - \bs{L} \leq 2 \vep$ and that $\Pr \li
\{ \bs{L} < p < \bs{U} \mid p \ri \}
> 1 - \de$ for any $p \in (0, 1)$ with prescribed $\vep \in (0, \f{1}{2})$ and $\de \in (0,
1)$.  For completeness, we shall develop multistage sampling schemes
in this setting.

Making use of the Clopper-Pearson confidence interval
\cite{Clopper}, we have established the following sampling scheme.
\beT \la{FW1} For $\al \in (0, 1)$ and integers $0 \leq k \leq n$,
define
\[ \mcal{L} (n, k, \al) = \left\{\begin{array}{ll}
   0 \;\;\;&  {\rm if}\; k=0\\
   \underline{p} \; \;\;\;&
   {\rm if}\; k > 0
\end{array} \right.\;\;\;{\rm and}\;\;\;
\mcal{U} (n, k, \al) = \left\{\begin{array}{ll}
   1  &  {\rm if}\; k= n\\
   \overline{p} \;
   & {\rm if}\; k < n
\end{array} \right.
\]
with $\underline{p} \in (0,1)$ satisfying $\sum_{j=k}^n {n \choose
j} \underline{p}^j (1- \underline{p})^{n-j} = \frac{\al}{2}$  and
$\overline{p} \in (0,1)$ satisfying $\sum_{j=0}^{k} {n \choose j}
   \overline{p}^j (1- \overline{p})^{n-j} = \frac{\al}{2}$.  Let $\ze > 0$ and $\ro > 0$. Let $n_1 < n_2 < \cd < n_s$ be the
ascending arrangement of all distinct elements of {\small $ \li \{
\li \lc \li ( \f{ 2 \vep^2} { \ln \f{1}{1- 2 \vep} } \ri )^{ 1 -
\f{i}{\tau} } \f{ \ln \f{1}{\ze \de} } { 2 \vep^2 } \ri \rc : i = 0,
1, \cd, \tau \ri \}$} with {\small  $\tau = \li \lc \f{ \ln \li (
\f{1}{ 2 \vep^2} \ln \f{1}{1- 2 \vep} \ri ) } { \ln (1 + \ro)} \ri
\rc$. }   For $\ell = 1, \cd, s$, define $K_\ell = \sum_{i =
1}^{n_\ell} X_i$  and $\bs{D}_\ell$ such that $\bs{D}_\ell = 1$ if
$\mcal{U} (n_\ell, K_\ell, \ze \de) - \mcal{L} (n_\ell, K_\ell, \ze
\de) \leq 2 \vep$; and $\bs{D}_\ell = 0$ otherwise.  Suppose the
stopping rule is that sampling is continued until $\bs{D}_\ell = 1$
for some $\ell \in \{1, \cd, s\}$.  Define {\small $\bs{L} =
\mcal{L} \li (\mathbf{n}, \sum_{i=1}^{\mathbf{n}} X_i, \ze \de \ri
)$} and {\small $\bs{U} = \mcal{U} \li (\mathbf{n},
\sum_{i=1}^{\mathbf{n}} X_i, \ze \de \ri )$}, where $\mathbf{n}$ is
the sample size when the sampling is terminated. Define {\small  \[
\mscr{Q}_L = \bigcup_{\ell = 1}^s \li \{ \mcal{L} (n_\ell, k, \ze
\de)  \in \li ( 0,  1 \ri ) : 0 \leq k \leq n_\ell \ri \}, \qqu
\mscr{Q}_U = \bigcup_{\ell = 1}^s \li \{ \mcal{U} (n_\ell, k, \ze
\de)  \in \li ( 0,  1 \ri ) : 0 \leq k \leq n_\ell \ri \}.
\]}
Then, a sufficient condition to guarantee $\Pr \li \{ \bs{L} < p <
\bs{U} \mid p \ri \}
> 1 - \de$ for any $p \in (0, 1)$ is that {\small  \bel &  &
\sum_{\ell = 1}^s \Pr \{ \mcal{L} (n_\ell, K_\ell, \ze \de) \geq p,
\; \bs{D}_{\ell - 1} = 0, \; \bs{D}_\ell = 1 \mid p \} < \f{\de}{2}
\qqu \fa p \in \mscr{Q}_L, \la{2D1F}\\
&  & \sum_{\ell = 1}^s \Pr \{ \mcal{U} (n_\ell, K_\ell, \ze \de)
\leq p, \; \bs{D}_{\ell - 1} = 0, \; \bs{D}_\ell = 1 \mid p \} <
\f{\de}{2} \qqu \fa p \in \mscr{Q}_U \la{2D2F}
 \eel}
where both (\ref{2D1F}) and (\ref{2D2F}) are satisfied if $0 < \ze <
\f{1}{2(\tau + 1)}$.

\eeT

Making use of Chernoff-Hoeffding inequalities \cite{Chernoff,
Hoeffding}, we have established the following sampling scheme.

\beT \la{FW2} For $\al \in (0, 1)$ and integers $0 \leq k \leq n$,
define
\[
\mcal{L} (n, k, \al) = \li \{\begin{array}{ll}
   \udl{p}  \;\;\;&  {\rm for}\; 0 < k < n,\\
   \li( \f{\al}{2} \ri)^{\f{1}{n}} \; \;\;\;&
   {\rm for} \; k = n,\\
   0 \; \;\;\;&
   {\rm for} \; k = 0
\end{array} \ri .  \qqu \mcal{U} (n, k, \al) = \li \{\begin{array}{ll}
   \ovl{p}  \;\;\;&  {\rm for}\; 0 < k < n,\\
   1 - \li( \f{\al}{2} \ri)^{\f{1}{n}} \; \;\;\;&
   {\rm for} \; k = 0,\\
   1 \; \;\;\;&
   {\rm for} \; k = n
\end{array} \ri .
\]
with $\udl{p} \in (0, \f{k}{n})$ satisfying $\mscr{M}_{\mrm{B}} \li
( \f{k}{n}, \udl{p} \ri )  = \f{ \ln (\ze \de) } { n}$ and $\ovl{p}
\in (\f{k}{n}, 1)$ satisfying $\mscr{M}_{\mrm{B}} \li ( \f{k}{n},
\ovl{p} \ri )  = \f{ \ln (\ze \de) } { n}$, where
$\mscr{M}_{\mrm{B}} (.,.)$ is a function such that
$\mscr{M}_{\mrm{B}} (z,\se) = z \ln \f{\se}{z} + (1 - z) \ln \f{1 -
\se}{1 - z}$ for $z \in (0,1)$ and $\se \in (0, 1)$. Let $\ze > 0$
and $\ro
> 0$. Let $n_1 < n_2 < \cd < n_s$ be the ascending arrangement of
all distinct elements of {\small $ \li \{ \li \lc \li ( \f{ 2
\vep^2} { \ln \f{1}{1- 2 \vep} } \ri )^{ 1 - \f{i}{\tau} } \f{ \ln
\f{1}{\ze \de} } { 2 \vep^2 } \ri \rc : i = 0, 1, \cd, \tau \ri \}$}
with {\small  $\tau = \li \lc \f{ \ln \li ( \f{1}{ 2 \vep^2} \ln
\f{1}{1- 2 \vep} \ri ) } { \ln (1 + \ro)} \ri \rc$. }   For $\ell =
1, \cd, s$, define $K_\ell = \sum_{i = 1}^{n_\ell} X_i$  and
$\bs{D}_\ell$ such that $\bs{D}_\ell = 1$ if $\mcal{U} (n_\ell,
K_\ell, \ze \de) - \mcal{L} (n_\ell, K_\ell, \ze \de) \leq 2 \vep$;
and $\bs{D}_\ell = 0$ otherwise.  Suppose the stopping rule is that
sampling is continued until $\bs{D}_\ell = 1$ for some $\ell \in
\{1, \cd, s\}$.  Define {\small $\bs{L} = \mcal{L} \li (\mathbf{n},
\sum_{i=1}^{\mathbf{n}} X_i, \ze \de \ri )$} and {\small $\bs{U} =
\mcal{U} \li (\mathbf{n}, \sum_{i=1}^{\mathbf{n}} X_i, \ze \de \ri
)$}, where $\mathbf{n}$ is the sample size when the sampling is
terminated. Define {\small  \[ \mscr{Q}_L = \bigcup_{\ell = 1}^s \li
\{ \mcal{L} (n_\ell, k, \ze \de)  \in \li ( 0,  1 \ri ) : 0 \leq k
\leq n_\ell \ri \}, \qqu \mscr{Q}_U = \bigcup_{\ell = 1}^s \li \{
\mcal{U} (n_\ell, k, \ze \de)  \in \li ( 0,  1 \ri ) : 0 \leq k \leq
n_\ell \ri \}.
\]}
Then, a sufficient condition to guarantee $\Pr \li \{ \bs{L} < p <
\bs{U} \mid p \ri \}
> 1 - \de$ for any $p \in (0, 1)$ is that {\small  \bel &  &
\sum_{\ell = 1}^s \Pr \{ \mcal{L} (n_\ell, K_\ell, \ze \de) \geq p,
\; \bs{D}_{\ell - 1} = 0, \; \bs{D}_\ell = 1 \mid p \} < \f{\de}{2}
\qqu \fa p \in \mscr{Q}_L, \la{2D1FC}\\
&  & \sum_{\ell = 1}^s \Pr \{ \mcal{U} (n_\ell, K_\ell, \ze \de)
\leq p, \; \bs{D}_{\ell - 1} = 0, \; \bs{D}_\ell = 1 \mid p \} <
\f{\de}{2} \qqu \fa p \in \mscr{Q}_U \la{2D2FC}
 \eel}
where both (\ref{2D1FC}) and (\ref{2D2FC}) are satisfied if $0 < \ze
< \f{1}{2(\tau + 1)}$.

\eeT

Making use of Massart's inequality \cite{Massart:90}, we have
established the following sampling scheme.

 \beT
\la{FW3} For $\al \in (0, 1)$ and integers $0 \leq k \leq n$, define
{\small \[  \mcal{L} (n, k, \al) = \max \li \{ 0, \; \frac{k}{n} +
\frac{3}{4} \; \frac{ 1 - \frac{2k}{n} - \sqrt{ 1 + \frac{9}{ 2 \ln
\frac{2}{\al} } \; k ( 1- \frac{k}{n}) } } {1 + \frac{9 n}{ 8 \ln
\frac{2}{\al} } } \ri \},
\]} and {\small \[ \mcal{U} (n, k, \al) = \min \li \{ 1, \; \frac{k}{n} +
\frac{3}{4} \; \frac{ 1 - \frac{2k}{n} + \sqrt{ 1 + \frac{9}{ 2 \ln
\frac{2}{\al} } \; k ( 1- \frac{k}{n}) } } {1 + \frac{9 n}{ 8 \ln
\frac{2}{\al} } } \ri \}.
\]}  Let $\ze > 0$ and $\ro
> 0$. Let $n_1 < n_2 < \cd < n_s$ be the ascending arrangement of
all distinct elements of {\small $ \li \{ \li \lc \f{8}{9} \li (
\frac{3}{4 \vep} + 1  \ri )^{ \f{i}{\tau} } \li ( \frac{3}{4 \vep} -
1  \ri ) \ln \frac{1}{\ze \de} \ri \rc : i = 0, 1, \cd, \tau \ri
\}$} with {\small  $\tau = \li \lc \f{ \ln \li ( \frac{3}{4 \vep} +
1 \ri ) } { \ln (1 + \ro)} \ri \rc$. } For $\ell = 1, \cd, s$,
define $K_\ell = \sum_{i = 1}^{n_\ell} X_i$ and $\bs{D}_\ell$ such
that $\bs{D}_\ell = 1$ if \[
1 - \frac{9}{ 2 \ln (\ze \de) } \;
K_\ell \li ( 1- \frac{K_\ell }{n_\ell} \ri ) \leq \vep^2 \li [
\f{4}{3} - \frac{3 n_\ell}{ 2 \ln (\ze \de) } \ri ]^2, \] and
$\bs{D}_\ell = 0$ otherwise.  Suppose the stopping rule is that
sampling is continued until $\bs{D}_\ell = 1$ for some $\ell \in
\{1, \cd, s\}$. Define {\small $\bs{L} = \mcal{L} \li (\mathbf{n},
\sum_{i=1}^{\mathbf{n}} X_i, \ze \de \ri )$} and {\small $\bs{U} =
\mcal{U} \li (\mathbf{n}, \sum_{i=1}^{\mathbf{n}} X_i, \ze \de \ri
)$}, where $\mathbf{n}$ is the sample size when the sampling is
terminated. Define {\small  \[ \mscr{Q}_L = \bigcup_{\ell = 1}^s \li
\{ \mcal{L} (n_\ell, k, \ze \de)  \in \li ( 0,  1 \ri ) : 0 \leq k
\leq n_\ell \ri \}, \qqu \mscr{Q}_U = \bigcup_{\ell = 1}^s \li \{
\mcal{U} (n_\ell, k, \ze \de)  \in \li ( 0,  1 \ri ) : 0 \leq k \leq
n_\ell \ri \}.
\]}
Then, a sufficient condition to guarantee $\Pr \li \{ \bs{L} < p <
\bs{U} \mid p \ri \}
> 1 - \de$ for any $p \in (0, 1)$ is that {\small  \bel &  &
\sum_{\ell = 1}^s \Pr \{ \mcal{L} (n_\ell, K_\ell, \ze \de) \geq p,
\; \bs{D}_{\ell - 1} = 0, \; \bs{D}_\ell = 1 \mid p \} < \f{\de}{2}
\qqu \fa p \in \mscr{Q}_L, \la{2D1mF}\\
&  & \sum_{\ell = 1}^s \Pr \{ \mcal{U} (n_\ell, K_\ell, \ze \de)
\leq p, \; \bs{D}_{\ell - 1} = 0, \; \bs{D}_\ell = 1 \mid p \} <
\f{\de}{2} \qqu \fa p \in \mscr{Q}_U \la{2D2mF}
 \eel}
where both (\ref{2D1mF}) and (\ref{2D2mF}) are satisfied if $0 < \ze
< \f{1}{2(\tau + 1)}$.

\eeT

\bsk

It should be noted that the interval estimation methods described in
Theorems \ref{FW1}--\ref{FW3} can be made less conservative by using
tight bounds of $C (p, \vep) = 1 - \Pr \{ \bs{L} < p < \bs{U}  \mid
p \}$ for $p \in [a, b] \subseteq \Se$ in Theorem \ref{bbsplit}.
Based on such bounds, a branch-and-bound type strategy described in
section 2.8 of \cite{Chen_EST} can be used to facilitate the search
of an appropriate value of $\ze$ such that the coverage probability
associated with interval $(\bs{L}, \bs{U})$ is no less than $1 -
\de$.

\beT  \la{bbsplit} Let $\mcal{L}_\ell = \mcal{L} (\wh{\bs{p}}_\ell,
\ze, \de)$ and $\mcal{U}_\ell = \mcal{U} (\wh{\bs{p}}_\ell, \ze,
\de)$ for $\ell = 1, \cd, s$.   Then, \bee C (p, \vep) & \leq & \Pr
\{  \bs{L} \geq a \mid b \} + \Pr \{  \bs{U} \leq b \mid a \} \\
& \leq & \sum_{\ell = 1}^s \Pr \{ \mcal{L}_\ell \geq a, \;
\bs{D}_{\ell - 1} = 0, \; \bs{D}_{\ell } = 1 \mid b \} + \sum_{\ell
= 1}^s \Pr \{ \mcal{U}_\ell \leq b, \; \bs{D}_{\ell - 1} = 0, \;
\bs{D}_{\ell } = 1 \mid a \}, \eee \bee C (p, \vep) & \geq & \Pr \{
\bs{L} \geq b \mid a \} + \Pr \{  \bs{U} \leq a \mid b \} \\
& \geq & \sum_{\ell = 1}^s \Pr \{ \mcal{L}_\ell \geq b, \;
\bs{D}_{\ell - 1} = 0, \; \bs{D}_{\ell } = 1 \mid a \} + \sum_{\ell
= 1}^s \Pr \{ \mcal{U}_\ell \leq a, \; \bs{D}_{\ell - 1} = 0, \;
\bs{D}_{\ell } = 1 \mid b \} \eee for any $p \in [a, b]$. Moreover,
if the open interval $(a, b)$ contains no element of the supports of
$\bs{L}$ and $\bs{U}$, then \bee C (p, \vep) & \leq & \Pr
\{  \bs{L} \geq b \mid b \} + \Pr \{  \bs{U} \leq a \mid a \} \\
& \leq & \sum_{\ell = 1}^s \Pr \{ \mcal{L}_\ell \geq b, \;
\bs{D}_{\ell - 1} = 0, \; \bs{D}_{\ell } = 1 \mid b \} + \sum_{\ell
= 1}^s \Pr \{ \mcal{U}_\ell \leq a, \; \bs{D}_{\ell - 1} = 0, \;
\bs{D}_{\ell } = 1 \mid a \}, \eee \bee C (p, \vep) & \geq & \Pr
\{  \bs{L} > a \mid a \} + \Pr \{  \bs{U} < b \mid b \} \\
& \geq  & \sum_{\ell = 1}^s \Pr \{ \mcal{L}_\ell > a, \;
\bs{D}_{\ell - 1} = 0, \; \bs{D}_{\ell } = 1 \mid a \} + \sum_{\ell
= 1}^s \Pr \{ \mcal{U}_\ell < b, \; \bs{D}_{\ell - 1} = 0, \;
\bs{D}_{\ell } = 1 \mid b \} \eee for any $p \in (a, b)$. \eeT

\bsk

We would like to note that Theorems 1 and 2 of \cite{Chen_EST} play
important roles in the establishment of the theorems in this
section. As can be seen from Theorems 1--6, the confidence
requirements can be satisfied by choosing $\ze$ to be sufficiently
small.  The application of the double-decision-variable method and
the single-decision-variable method is obvious.  To determine $\ze$
as large as possible and thus make the sampling schemes most
efficient, the computational techniques such as bisection confidence
tuning, domain truncation, triangular partition developed in
\cite{Chen_EST} can be applied.

With regard to the tightness of the double-decision-variable method,
we can develop results similar to Theorems 13, 18 and 23 of
\cite{Chen_EST}.

With regard to the asymptotic performance of our sampling schemes,
we can develop results similar to Theorems 14, 19 and 24 of
\cite{Chen_EST}.

\sect{Estimation of Bounded-variable Means}

 The method proposed for estimating binomial parameters can be
 generalized for estimating means of random variables bounded in interval $[0, 1]$.
 Formally, let $Z \in [0, 1]$ be a random variable with expectation $\mu =
 \bb{E} [Z]$. We can estimate $\mu$ based on i.i.d. random samples
 $Z_1, Z_2, \cd$ of $Z$ by virtue of the following results.

 \beT \la{coverage_abs} Let $0 < \vep < \f{1}{2}$ and $0 < \de < 1$.
Let $n_1 < n_2 < \cd < n_s$ be a sequence of sample sizes such that
{\small $n_s \geq \f{ \ln \f{2 s} { \de} } { 2 \vep^2 }$}.  Define
$\wh{\bs{\mu}}_\ell = \f{ \sum_{i=1}^{n_\ell} Z_i }{n_\ell}$  for
$\ell = 1, \cd, s$.  Suppose the stopping rule is that sampling is
continued until {\small $\li ( \li | \wh{\bs{\mu}}_\ell - \f{1}{2}
\ri | - \f{2 \vep }{3} \ri )^2 \geq \f{1}{4} - \f{ \vep^2 n_\ell }
{2 \ln (2 s \sh \de) }$} for some $\ell \in \{1, \cd, s \}$. Define
$\bs{\wh{\mu}} = \f{\sum_{i=1}^{\mathbf{n}} Z_i}{\mathbf{n}}$ where
$\mathbf{n}$ is the sample size when the sampling is terminated.
 Then, $\Pr \li \{ \li | \bs{\wh{\mu}} - \mu \ri | < \vep \ri \}
\geq 1 - \de$. \eeT

\bsk

This theorem can be shown by a variation of the argument for Theorem
1.

\beT \la{coverage_mix_bound} Let $0 < \de < 1, \; 0 < \vep_a <
\f{3}{8}$ and {\small $\f{6 \vep_a}{3 - 2 \vep_a } < \vep_r < 1$}.
Let $n_1 < n_2 < \cd < n_s$ be a sequence of sample sizes such that
{\small $n_s \geq 2  \li ( \f{1}{ \vep_r } + \frac{1} {3} \ri ) \li
( \f{1}{ \vep_a } - \f{1}{ \vep_r } - \frac{1} {3} \ri ) \ln \li (
\f{2s}{\de} \ri )$}.  Define $\wh{\bs{\mu}}_\ell = \f{
\sum_{i=1}^{n_\ell} Z_i }{n_\ell}$  for $\ell = 1, \cd, s$.  Define
{\small \[ \bs{D}_\ell = \bec 0 & \mrm{for} \; \f{1}{2} - \f{2}{3}
\vep_a - \sq{ \f{1}{4} + \f{ n_\ell \vep_a^2 } {2 \ln (\ze \de) } }
< \wh{\bs{\mu}}_\ell < \f{ 6(1 - \vep_r) (3 - \vep_r) \ln (\ze
\de) } { 2 (3 - \vep_r)^2 \ln (\ze \de) - 9 n_\ell \vep_r^2} \; \mrm{or}\\
  &  \qu \;\; \f{1}{2} + \f{2}{3} \vep_a - \sq{ \f{1}{4} + \f{
n_\ell \vep_a^2 } {2 \ln (\ze \de) } } < \wh{\bs{\mu}}_\ell < \f{
6(1 + \vep_r) (3 + \vep_r) \ln (\ze \de) } { 2 (3 +
\vep_r)^2 \ln (\ze \de) - 9 n_\ell  \vep_r^2},\\
1 & \mrm{else}
 \eec
\]}
for $\ell = 1, \cd, s - 1$ and $\bs{D}_s = 1$.  Suppose the stopping
rule is that sampling is continued until $\bs{D}_\ell = 1$ for some
$\ell \in \{1, \cd, s \}$.  Define $\bs{\wh{\mu}} =
\f{\sum_{i=1}^{\mathbf{n}} Z_i}{\mathbf{n}}$ where $\mathbf{n}$ is
the sample size when the sampling is terminated.
 Then, $\Pr \li \{ \li | \bs{\wh{\mu}} - \mu \ri | < \vep_a \; \tx{or} \; \li | \bs{\wh{\mu}} - \mu \ri | < \vep_r \mu \ri \}
\geq 1 - \de$. \eeT

\bsk

This theorem can be shown by a variation of the argument for Theorem
2.  In the general case that $Z$ is a random variable bounded in
$[a, b]$, it is useful to estimate the mean $\mu = \bb{E} [ Z]$
based on i.i.d. samples of $Z$ with a mixed criterion.  For this
purpose,  we shall introduce the function
\[
\mcal{M}(z, \mu) = \bec  \f{ (\mu - z)^2 } {2 \li ( \f{2 \mu}{3} +
\f{z}{3} \ri ) \li ( \f{2 \mu}{3} + \f{z}{3} - 1 \ri ) } & \tx{for}
\; 0 \leq z \leq 1 \; \tx{and} \;  \mu \in (0, 1),\\
- \iy & \tx{for} \; 0 \leq z \leq 1 \; \tx{and} \; \mu \notin (0, 1)
\eec
\]
and propose the following multistage estimation method.

\beT \la{coverage_mix_general} Let $0 < \de < 1, \; \vep_a > 0$ and
$0 < \vep_r < 1$. Let $n_1 < n_2 < \cd < n_s$ be a sequence of
sample sizes such that {\small $n_s \geq  \f{(b - a)^2}{ 2 \vep_a^2}
\ln \li ( \f{2s}{\de} \ri )$}. Define $\wh{\bs{\mu}}_\ell = \f{
\sum_{i=1}^{n_\ell} Z_i }{n_\ell}, \; \wt{\bs{\mu}}_\ell  = a +
\f{1}{b - a} \wh{\bs{\mu}}_\ell$,
\[
 \udl{\bs{\mu}}_\ell = a +
\f{1}{b - a}  \min \li \{ \wh{\bs{\mu}}_\ell - \vep_a, \;
\f{\wh{\bs{\mu}}_\ell}{ 1 + \mrm{sgn} (\wh{\bs{\mu}}_\ell) \vep_r }
\ri \}, \qqu \ovl{\bs{\mu}}_\ell = a + \f{1}{b - a}  \max \li \{
\wh{\bs{\mu}}_\ell + \vep_a, \; \f{\wh{\bs{\mu}}_\ell}{ 1 -
\mrm{sgn} (\wh{\bs{\mu}}_\ell) \vep_r } \ri \}
\]
for $\ell = 1, \cd, s$.  Suppose the stopping rule is that sampling
is continued until $\mcal{M}(\wt{\bs{\mu}}_\ell,
\udl{\bs{\mu}}_\ell) \leq \f{1}{n_\ell} \ln \f{\de}{2s}$  and
$\mcal{M}(\wt{\bs{\mu}}_\ell, \ovl{\bs{\mu}}_\ell) \leq
\f{1}{n_\ell} \ln \f{\de}{2s}$ for some $\ell \in \{1, \cd, s \}$.
Define $\bs{\wh{\mu}} = \f{\sum_{i=1}^{\mathbf{n}} Z_i}{\mathbf{n}}$
where $\mathbf{n}$ is the sample size when the sampling is
terminated.
 Then, $\Pr \li \{ \li | \bs{\wh{\mu}} - \mu \ri | < \vep_a \; \tx{or} \; \li | \bs{\wh{\mu}} - \mu \ri | < \vep_r |\mu| \ri \}
\geq 1 - \de$. \eeT

 \sect{A Link between
Binomial and Bounded Variables}

There exists an inherent connection between a binomial parameter and
the mean of a bounded variable.  In this regard, we have

\beT \la{Link} Let $Z$ be a random variable bounded in $[0, 1]$. Let
$U$ a random variable uniformly distributed over $[0, 1]$. Suppose
$Z$ and $U$ are independent. Then,
\[ \bb{E} [Z] = \Pr \{ Z \geq U \}.
\] \eeT

\bpf Let $F_{Z, U}$ be the joint distribution of $Z$ and $U$.  Let
$F_Z$ be the cumulative distribution function of $Z$. Since $Z$ and
$U$ are independent, using Riemann-Stieltjes integration, we have
\[
\Pr \{ Z \geq U \} = \int_{z = 0}^1 \int_{u = 0}^z d F_{Z, U} =
\int_{z = 0}^1 \int_{u = 0}^z d u \; d F_Z =  \int_{z = 0}^1 z \; d
F_Z = \bb{E} [Z].
\]
\epf

To see why Theorem \ref{Link} reveals a relationship between the
mean of a bounded variable and a binomial parameter, we define
\[
X = \bec 1 & \tx{for} \; Z \geq U,\\
0 & \tx{otherwise}. \eec
\]
Then, by Theorem \ref{Link}, we have $\Pr \{ X = 1 \} = 1 - \Pr \{ X
= 0 \} = \bb{E} [Z]$. This implies that $X$ is a Bernoulli random
variable and $\bb{E} [Z]$ is actually a binomial parameter. As a
consequence, the techniques of estimating a binomial parameter can
be useful for estimating the mean of a bounded variable. Specially,
for a sequence of i.i.d. random samples $Z_1, Z_2, \cd$ of bounded
variable $Z$ and a sequence of i.i.d. random samples $U_1, U_2, \cd$
of uniform variable $U$ such that that $Z_i$ is independent with
$U_i$ for all $i$, we can define a sequence of i.i.d. random samples
$X_1, X_2, \cd$  of Bernoulli random variable $X$ by
\[
X_i = \bec 1 & \tx{for} \; Z_i \geq U_i,\\
0 & \tx{otherwise}. \eec
\]

\sect{Conclusion}

We have established a new multistage approach for estimating the
mean of a bounded variable.  Our approach can provide an estimator
for the unknown mean which rigorously guarantees prescribed levels
of precision and confidence. Our approach is also very flexible in
the sense that the precision can be expressed in terms of different
types of margins of errors.

\bsk

\appendix

\sect{Preliminary Results for Proofs of Theorems}

 We need some preliminary results, especially some properties of function $\mcal{M}(z, \mu)$ defined in Section 3.

\beL \la{decrea} $\mcal{M} (z, z + \vep )$ is monotonically
increasing with respect to  $z \in (0, \f{1}{2} - \f{2 \vep}{3})$,
and is monotonically decreasing with respect to $z \in (\f{1}{2} -
\f{2 \vep}{3}, 1 - \vep)$.  Similarly, $\mcal{M} (z, z - \vep )$ is
monotonically increasing with respect to $z \in (\vep, \f{1}{2} +
\f{2 \vep}{3})$, and is monotonically decreasing with respect to $z
\in ( \f{1}{2} + \f{2 \vep}{3}, 1)$. \eeL

\bpf The lemma can be established by checking the partial
derivatives
\[
\f{\pa \mcal{M} (z, z + \vep )} { \pa z } = \f{ \vep^2 } { \li [ \li
(z + \f{2 \vep}{3} \ri ) \li ( 1- z - \f{2 \vep}{3} \ri ) \ri ]^2 }
\li ( \f{1}{2} - \f{2 \vep}{3} - z \ri ),
\]
\[
\f{\pa \mcal{M} (z, z - \vep )} { \pa z } = \f{ \vep^2 } { \li [ \li
(z - \f{2 \vep}{3} \ri ) \li ( 1- z + \f{2 \vep}{3} \ri ) \ri ]^2 }
\li ( \f{1}{2} + \f{2 \vep}{3} - z \ri ).
\]

\epf

\beL \la{lem888m}  Let $0 < \vep < \f{1}{2}$. Then, $\mcal{M}(z, z +
\vep) \geq  \mcal{M}(z, z - \vep)$ for $z \in \li [0, \f{1}{2} \ri
]$, and $\mcal{M}(z, z + \vep) < \mcal{M}(z, z - \vep)$ for $z \in
\li ( \f{1}{2}, 1 \ri ]$. \eeL

\bpf

By the definition of the function $\mcal{M} (.,.)$, we have that
$\mcal{M} (z, \mu) = - \iy$ for $z \in [0, 1]$ and $\mu \notin (0,
1)$. Hence, the lemma is trivially true for $0 \leq z \leq \vep$ or
$1 - \vep \leq z \leq 1$.  It remains to show the lemma for $z \in
(\vep, 1 - \vep)$.  This can be accomplished by noting that
\[
\mcal{M} (z, z + \vep ) - \mcal{M} (z, z - \vep ) = \f{ 2 \vep^3 ( 1
- 2 z) } {3 \li (z + \f{2 \vep}{3} \ri ) \li ( 1- z - \f{2 \vep}{3}
\ri )  \li (z - \f{2 \vep}{3} \ri ) \li ( 1- z + \f{2 \vep}{3} \ri )
}.
\]
where the right-hand side is seen to be positive for $z \in \li (
\vep, \f{1}{2} \ri )$ and negative for $z \in \li ( \f{1}{2}, 1 -
\vep \ri )$. \epf

\beL \la{lem888m2B}

$\mcal{M} \li ( z, \f{z}{1 + \vep} \ri ) > \mcal{M} \li ( z, \f{z}{1
- \vep} \ri )$ for $0 < z < 1 - \vep < 1$.

\eeL

\bpf

It can be verified that
\[
\mcal{M} \li ( z, \f{z}{1 + \vep} \ri ) - \mcal{M} \li ( z, \f{z}{1
- \vep} \ri ) = \f{ 2 \vep^3 z (2 - z) } {3 \li ( 1 + \f{\vep}{3}
\ri ) \li [ 1 - z + \vep \li ( 1 - \f{ z}{3} \ri ) \ri ]  \li ( 1 -
\f{\vep}{3} \ri ) \li [ 1 - z - \vep \li ( 1 - \f{ z}{3} \ri ) \ri
]},
\]
from which it can be seen that {\small $\mcal{M} \li ( z, \f{z}{1 +
\vep} \ri ) > \mcal{M} \li ( z, \f{z}{1 - \vep} \ri )$} for $z \in
(0, 1 - \vep)$.

\epf

\beL \la{lemm22m}

$\mcal{M} (\mu - \vep, \mu) < \mcal{M} (\mu + \vep, \mu)$ for $0 <
\vep < \mu < \f{1}{2} < 1 - \vep$.

\eeL

\bpf The lemma follows from the fact that
\[
\mcal{M} (\mu - \vep, \mu) - \mcal{M} (\mu + \vep, \mu) = \f{ \vep^3
(2 \mu - 1) } { 3 \li ( \mu - \f{\vep}{3} \ri ) \li ( 1 - \mu +
\f{\vep}{3} \ri )  \li ( \mu + \f{\vep}{3} \ri ) \li ( 1 - \mu -
\f{\vep}{3} \ri )  },
\]
where the right-hand side is negative for $0 < \vep < \mu < \f{1}{2}
< 1 - \vep$.

\epf

\beL \la{decrev}
{\small $\mcal{M} \li ( z, \f{z}{1 + \vep} \ri )$}
is monotonically decreasing with respect to $z \in (0, 1)$.
Similarly, {\small $\mcal{M} \li ( z, \f{z}{1 - \vep} \ri )$} is
monotonically decreasing with respect to $z \in (0, 1 - \vep)$. \eeL

\bpf

The lemma can be shown by verifying that
\[
\f{\pa } { \pa z} \mcal{M} \li ( z, \f{z}{1 + \vep} \ri ) = - \f{
\vep^2}{2 \li ( 1 + \f{\vep}{3} \ri )} \times \f{ 1 + \vep } { \li [
(1 + \vep) ( 1 - z) + \f{2 \vep z}{3} \ri ]^2 } < 0
\]
for $z \in (0, 1)$ and that
\[
\f{\pa } { \pa z} \mcal{M} \li ( z, \f{z}{1 - \vep} \ri ) = - \f{
\vep^2}{2 \li ( 1 - \f{\vep}{3} \ri ) } \times \f{ 1 - \vep } { \li
[ (1 - \vep) ( 1 - z) - \f{2 \vep z}{3} \ri ]^2 } < 0
\]
for $z \in (0, 1 - \vep)$.

\epf

\beL

\la{Ddec} For any fixed $z \in (0, 1)$, $\mcal{M}(z,\mu)$ is
monotonically increasing with respect to $\mu \in (0, z)$, and is
monotonically decreasing with respect to $\mu \in (z, 1)$.
Similarly, for any fixed $\mu \in (0, 1)$, $\mcal{M}(z,\mu)$ is
monotonically increasing with respect to $z \in (0, \mu)$, and is
monotonically decreasing with respect to $z \in (\mu, 1)$.

\eeL

\bpf

The lemma can be shown by checking the following partial
derivatives:
\bee \f{\pa \mcal{M}(z,\mu)}{\pa \mu} & = & \f{(z -
\mu) \li [ \mu(1 - z) + z ( 1 - \mu) + z ( 1 -z) \ri ] } { 3 \li [
\li ( \f{2\mu}{3} + \f{z}{3} \ri ) \li ( 1 - \f{2\mu}{3} - \f{z}{3}
\ri ) \ri ]^2
},\\
\f{\pa \mcal{M}(z,\mu)}{\pa z} & = & \f{(\mu - z) \li [ \mu(1 -
\f{2\mu}{3} - \f{z}{3} ) + \f{z - \mu}{6} \ri ] } {
 \li [ \li ( \f{2\mu}{3} + \f{z}{3} \ri ) \li ( 1 - \f{2\mu}{3} - \f{z}{3} \ri ) \ri ]^2 }
  = \f{(\mu - z) \li [ (1 - \mu) (\f{2\mu}{3} + \f{z}{3} ) + \f{\mu - z}{6} \ri ] } {
 \li [ \li ( \f{2\mu}{3} + \f{z}{3} \ri ) \li ( 1 - \f{2\mu}{3} - \f{z}{3} \ri ) \ri ]^2 }.
\eee

\epf

The following result, stated as Lemma \ref{lem12},  is due to
Massart \cite{Massart:90}.

\beL \la{lem12} Let $\ovl{X}_n = \f{\sum_{i=1}^n X_i}{n}$ where
$X_1, \cd, X_n$ are i.i.d. random variables such that $0 \leq X_i
\leq 1$ and $\bb{E}[X_i] = \mu \in (0,1)$ for $i = 1, \cd, n$. Then,
 {\small ${\Pr} \li \{ \ovl{X}_n \geq z
\ri \}  <  \exp \li ( n \mcal{M} (z, \mu) \ri )$ for any $z \in
(\mu, 1)$. Similarly, ${\Pr} \li \{ \ovl{X}_n \leq z \ri \} < \exp
\li ( n \mcal{M} (z, \mu)  \ri )$} for any $z \in (0, \mu)$. \eeL

\beL \la{lem1} Let $\ovl{X}_n = \f{\sum_{i=1}^n X_i}{n}$ where $X_1,
\cd, X_n$ are i.i.d. random variables such that $0 \leq X_i \leq 1$
and $\bb{E}[X_i] = \mu \in (0,1)$ for $i = 1, \cd, n$. Then, {\small
$\Pr \li \{ \ovl{X}_n \geq \mu, \; \mcal{M} \li ( \ovl{X}_n ,  \mu
\ri ) \leq \f{\ln \al}{n} \ri \} \leq \al$} for any $\al > 0$. \eeL

\bpf  Since the lemma is trivially true for $\al \geq 1$, it remains
to show it for $\al \in (0, 1)$.  It can be checked that $\mcal{M}
(\mu, \mu ) = 1$ and {\small $\mcal{M} (1, \mu ) = \f{ 9 (\mu - 1) }
{ 4 (2 \mu + 1) }$}. Since {\small $\f{\pa \mcal{M} (z, \mu) } { \pa
z} = (\mu - z) [ \mu(1 - \f{2\mu}{3} - \f{z}{3} ) + \f{z - \mu}{6} ]
\sh
 [ ( \f{2\mu}{3} + \f{z}{3} ) ( 1 - \f{2\mu}{3} - \f{z}{3} ) ]^2   < 0$} for $z \in (\mu,1)$,
 we have that $\mcal{M} (z, \mu)$ is monotonically decreasing from
 $0$ to {\small $\f{ 9 (\mu - 1) } { 4 (2 \mu + 1) }$} as $z$ increases from $\mu$ to $1$.
 To show the lemma, we need to consider three cases as follows.

Case (i): {\small $ \f{ 9 (\mu - 1) } { 4 (2 \mu + 1) } > \f{\ln
\al}{n}$}. In this case, we have that {\small $\li \{ \ovl{X}_n \geq
\mu , \;  \mcal{M} \li ( \ovl{X}_n ,  \mu \ri ) \leq \f{\ln \al}{n}
\ri \}$} is an impossible event and the corresponding probability is
$0$.  This is because the minimum of $\mcal{M} (z, \mu)$ with
respect to $z \in (\mu, 1]$ is equal to $\f{ 9 (\mu - 1) } { 4 (2
\mu + 1) }$, which is greater than $\f{\ln \al}{n}$.

Case (ii): {\small $ \f{ 9 (\mu - 1) } { 4 (2 \mu + 1) } = \f{\ln
\al}{n}$}.  In this case, we have that {\small $\Pr \li \{ \ovl{X}_n
\geq \mu , \;  \mcal{M} \li ( \ovl{X}_n ,  \mu \ri ) \leq \f{\ln
\al}{n} \ri \} =  \Pr \{ X_i = 1, \; i = 1, \cd, n \} = \prod_{i =
1}^n \Pr \{ X_i = 1 \} \leq \prod_{i = 1}^n \bb{E} [ X_i ] = \mu^n <
\exp \li ( n  \times \f{ 9 (\mu - 1) } { 4 (2 \mu + 1) } \ri ) =
\al$}, where the last inequality is due to the fact that {\small
$\ln \mu < \f{ 9 (\mu - 1) } { 4 (2 \mu + 1) }$}.  To prove this
fact,  we define $g(\mu) = \ln \mu - \f{ 9 (\mu - 1) } { 4 (2 \mu +
1) }$. Then, the first derivative of $g(\mu)$ is $g^\prime (\mu) =
\f{ 5 \mu^2 + 4 - 11 \mu (1 - \mu)  } { 4 \mu \li ( 2 \mu + 1 \ri
)^2 } \geq \f{ 5 \mu^2 + 4 - 11 \times \f{1}{4} } { 4 \mu \li ( 2
\mu + 1 \ri )^2 } > 0$ for any $\mu \in (0, 1)$. This implies that
$g(\mu)$ is monotonically increasing with respect to $\mu \in (0,
1)$. By virtue of such monotonicity and the fact that $g(1) = 0$, we
can conclude that $g (\mu) < 0$ for any $\mu \in (0, 1)$.  This
establishes {\small $\ln \mu < \f{ 9 (\mu - 1) } { 4 (2 \mu + 1)
}$}.

Case (iii): {\small $ \f{ 9 (\mu - 1) } { 4 (2 \mu + 1) }  < \f{\ln
\al}{n}$}. In this case, there exists a unique number $z^* \in (\mu
, 1)$ such that $\mcal{M} (z^*, \mu) = \f{\ln \al}{n}$. Since
$\mcal{M} (z, \mu )$ is monotonically decreasing with respect to $z
\in (\mu, 1)$, it must be true that any $\ovl{x} \in (\mu, 1)$
satisfying $\mcal{M} (\ovl{x}, \mu ) \leq \f{\ln \al}{n}$ is no less
than $z^*$.  This implies that {\small $\li \{ \ovl{X}_n \geq \mu,
\; \mcal{M} \li ( \ovl{X}_n, \mu \ri ) \leq \f{\ln \al}{n} \ri \}
\subseteq \li \{ \ovl{X}_n \geq z^* \ri \}$} and {\small $\Pr \li \{
\ovl{X}_n \geq \mu , \; \mcal{M} \li ( \ovl{X}_n, \mu \ri ) \leq
\f{\ln \al}{n} \ri \}
 \leq  \Pr \li \{ \ovl{X}_n \geq z^* \ri \}  \leq  \exp (n
\; \mcal{M} (z^*, \mu)) = \al$}, where the last inequality follows
from Lemma \ref{lem12}. This completes the proof of the lemma.

\epf

\beL \la{lem2} Let $\ovl{X}_n = \f{\sum_{i=1}^n X_i}{n}$ where $X_1,
\cd, X_n$ are i.i.d. random variables such that $0 \leq X_i \leq 1$
and $\bb{E}[X_i] = \mu \in (0,1)$ for $i = 1, \cd, n$. Then, {\small
$\Pr \{ \ovl{X}_n \leq \mu, \;  \mcal{M} \li ( \ovl{X}_n ,  \mu \ri
) \leq \f{\ln \al}{n} \}  \leq \al$} for any $\al > 0$. \eeL

\bpf Since the lemma is trivially true for $\al \geq 1$, it suffices
to show it for $\al \in (0, 1)$.  It can be checked that $\mcal{M} (
\mu, \mu) = 1$ and {\small $\mcal{M} (0, \mu) = \f{ 9 \mu  } { 4 (2
\mu - 3) }$}. Since {\small $\f{\pa \mcal{M} (z, \mu) } { \pa z}  =
(\mu - z) [ (1 - \mu) (\f{2\mu}{3} + \f{z}{3} ) + \f{\mu - z}{6} ]
\sh [ ( \f{2\mu}{3} + \f{z}{3} ) ( 1 - \f{2\mu}{3} - \f{z}{3}  ) ]^2
> 0$} for $z \in (0, \mu)$, we have that $\mcal{M} (z, \mu)$ is
monotonically increasing from $\f{ 9 \mu  } { 4 (2 \mu - 3) }$ to
$0$ as $z$ increases from $0$ to $\mu$.  Now there are three cases:

Case (i): {\small $ \f{ 9 \mu  } { 4 (2 \mu - 3) } > \f{\ln
\al}{n}$}. In this case, we have that {\small $\li \{ \ovl{X}_n \leq
\mu , \; \mcal{M} \li ( \ovl{X}_n ,  \mu \ri ) \leq \f{\ln \al}{n}
\ri \}$} is an impossible event and the corresponding probability is
$0$.  This is because the minimum of $\mcal{M} (z, \mu)$ with
respect to $z \in [0, \mu )$ is equal to $\f{ 9 \mu} { 4 (2 \mu - 3)
}$, which is greater than $\f{\ln \al}{n}$.

Case (ii): {\small $ \f{ 9 \mu  } { 4 (2 \mu - 3) } = \f{\ln
\al}{n}$}.  In this case, we have that {\small $\Pr \li \{ \ovl{X}_n
\leq \mu, \; \mcal{M} \li ( \ovl{X}_n ,  \mu \ri ) \leq \f{\ln
\al}{n} \ri \} = \Pr \{ X_i = 0, \; i = 1, \cd, n \} = \prod_{i =
1}^n \Pr \{ X_i = 0 \} = \prod_{i = 1}^n (1 - \Pr \{ X_i \neq 0 \} )
\leq \prod_{i = 1}^n (1 - \bb{E} [X_i] ) = (1 - \mu)^n < \exp \li (
n \times \f{ 9 \mu  } { 4 (2 \mu - 3) } \ri ) = \al $}, where the
last inequality is due to the fact that {\small $\ln (1 - \mu) < \f{
9 \mu } { 4 (2 \mu - 3) }$}.  To prove this fact,  we define $h(\mu)
= \ln (1 - \mu) - \f{ 9 \mu } { 4 (2 \mu - 3) }$. Then, the first
derivative of $h(\mu)$ is $h^\prime (\mu) = \f{ - 16 \mu^2 + 21 \mu
- 9  } { 4 (1 - \mu) \li ( 2 \mu - 3 \ri )^2 } \leq \f{ 16 \times (
\f{21}{32})^2 - 9  } { 4 (1 - \mu) \li ( 2 \mu - 3 \ri )^2 } < 0$
for any $\mu \in (0, 1)$. This implies that $h(\mu)$ is
monotonically decreasing with respect to $\mu \in (0, 1)$. By virtue
of such monotonicity and the fact that $h(0) = 0$, we can conclude
that $h (\mu) < 0$ for any $\mu \in (0, 1)$. This establishes
{\small $\ln (1 - \mu) < \f{ 9 \mu } { 4 (2 \mu - 3) }$}.

Case (iii): {\small $ \f{ 9 \mu  } { 4 (2 \mu - 3) } < \f{\ln
\al}{n}$}.  In this case, there exists a unique number $z^* \in (0,
\mu)$ such that $\mcal{M} (z^*, \mu) = \f{\ln \al}{n}$. Since
$\mcal{M} (z, \mu)$ is monotonically increasing with respect to $z
\in (0, \mu)$, it must be true that any $\ovl{x} \in (0, \mu)$
satisfying $\mcal{M} (\ovl{x}, \mu) \leq \f{\ln \al}{n}$ is no
greater than $z^*$.  This implies that {\small $\li \{ \ovl{X}_n
\leq  \mu, \; \mcal{M} \li ( \ovl{X}_n ,  \mu \ri ) \leq \f{\ln
\al}{n} \ri \} \subseteq \li \{ \ovl{X}_n \leq z^* \ri \}$} and thus
{\small $\Pr \li \{ \ovl{X}_n \leq \mu, \; \mcal{M} \li ( \ovl{X}_n
, \mu \ri ) \leq \f{\ln \al}{n} \ri \} \leq \Pr \li \{ \ovl{X}_n
\leq z^* \ri \} \leq \exp ( n \mcal{M} ( z^*, \mu) ) = \al$ }, where
the last inequality follows from Lemma \ref{lem12}.  This completes
the proof of the lemma.

\epf

\sect{Proof of Theorem 1}

Throughout the proof of Theorem 1, we define random variables
$\bs{D}_\ell, \; \ell = 1, \cd, s$ such that $\bs{D}_\ell = 1$ if
{\small $\li ( \li | \wh{\bs{p}}_\ell - \f{1}{2} \ri | -
\f{2\vep}{3} \ri )^2 \geq \f{1}{4} + \f{  n_\ell \; \vep^2 } {2 \ln
(\ze \de) }$} and $\bs{D}_\ell = 0$ otherwise.  Then, the stopping
rule can be restated as ``sampling is continued until $\bs{D}_\ell =
1$ for some $\ell \in \{1, \cd, s \}$''.

\bsk

\beL \la{absD}

$\bs{D}_s = 1$.

\eeL \bpf

By the definition of $\bs{D}_s$, we have that {\small $\{ \bs{D}_s =
1 \} = \li \{ \li ( \li | \wh{\bs{p}}_s - \f{1}{2} \ri | -
\f{2\vep}{3} \ri )^2 \geq \f{1}{4} + \f{  n_s \; \vep^2 } {2 \ln
(\ze \de) } \ri \}$}.  By the definition of sample sizes, we have
{\small $n_s = \li \lc  \f{ \ln \f{1}{\ze \de} } { 2 \vep^2  } \ri
\rc \geq \f{ \ln \f{1}{\ze \de} } { 2 \vep^2  }$}, which implies
that $\f{1}{4} + \f{  n_s \; \vep^2 } {2 \ln (\ze \de) } \leq 0$.
Since $\{ \li ( \li | \wh{\bs{p}}_s - \f{1}{2} \ri | - \f{2\vep}{3}
\ri )^2 \geq 0 \}$ is a sure event, it follows that {\small $\li \{
\li ( \li | \wh{\bs{p}}_s - \f{1}{2} \ri | - \f{2\vep}{3} \ri )^2
\geq \f{1}{4} + \f{  n_s \; \vep^2 } {2 \ln (\ze \de) } \ri \}$} is
a sure event and consequently $\bs{D}_s = 1 $. This completes the
proof of the lemma.

\epf

\beL \la{abs381}
{\small $\li \{  \wh{\bs{p}}_\ell \leq
  p - \vep, \; \bs{D}_\ell = 1 \ri \}  \leu \li \{ \wh{\bs{p}}_\ell < p,  \; \mcal{M} \li ( \wh{\bs{p}}_\ell, p \ri )
\leq \f{\ln (\ze \de)}{n_\ell} \ri \} $} for $\ell = 1, \cd, s$.

\eeL

\bpf

Since {\small $\{ \bs{D}_\ell = 1 \} = \li \{ \li ( \li |
\wh{\bs{p}}_\ell - \f{1}{2} \ri | - \f{2\vep}{3} \ri )^2 \geq
 \f{1}{4} + \f{  n_\ell \vep^2 } {2 \ln (\ze \de) } \ri \}$}, it suffices to show {\small
\[
\li \{
\wh{\bs{p}}_\ell \leq p - \vep,  \; \li ( \li | \wh{\bs{p}}_\ell -
\f{1}{2} \ri | - \f{2\vep}{3} \ri )^2 \geq
 \f{1}{4} + \f{  n_\ell \vep^2 } {2 \ln (\ze \de) } \ri \} \subseteq \li \{ \wh{\bs{p}}_\ell < p, \;
\mcal{M} \li ( \wh{\bs{p}}_\ell, p \ri ) \leq \f{\ln (\ze
\de)}{n_\ell} \ri \}. \]
 }
For this purpose, we let {\small $\om \in \li \{ \wh{\bs{p}}_\ell
\leq p - \vep, \; \li ( \li | \wh{\bs{p}}_\ell - \f{1}{2} \ri | -
\f{2\vep}{3} \ri )^2 \geq \f{1}{4} + \f{  n_\ell \vep^2 } {2 \ln
(\ze \de) } \ri \}, \; \wh{p}_\ell = \wh{\bs{p}}_\ell (\om)$} and
proceed to show $\wh{p}_\ell < p, \; \mcal{M} \li ( \wh{p}_\ell, p
\ri ) \leq \f{\ln (\ze \de)}{n_\ell}$. Clearly, $\wh{p}_\ell < p$
follows immediately from $\wh{p}_\ell \leq p - \vep$.  To show
$\mcal{M} \li ( \wh{p}_\ell, p \ri ) \leq \f{\ln (\ze
\de)}{n_\ell}$, we need to establish \be \la{eqDD} \li ( \wh{p}_\ell
- \f{1}{2}  + \f{2\vep}{3} \ri )^2 \geq \f{1}{4} + \f{  n_\ell
\vep^2 } {2 \ln (\ze \de) } \ee based on \be \la{eqD}
 \li ( \li |
\wh{p}_\ell - \f{1}{2} \ri | - \f{2\vep}{3} \ri )^2 \geq
 \f{1}{4} + \f{  n_\ell \vep^2 } {2 \ln (\ze \de) }. \ee  It is obvious that (\ref{eqDD}) holds if
{\small $\f{1}{4} + \f{  n_\ell \vep^2 } {2 \ln (\ze \de) } \leq
0$}.  It remains to show (\ref{eqDD}) under the condition that
{\small $\f{1}{4} + \f{  n_\ell \vep^2 } {2 \ln (\ze \de) }
> 0$}.  Note that (\ref{eqD}) implies either \be \la{eqD1} \li |
\wh{p}_\ell - \f{1}{2} \ri | - \f{2\vep}{3} \geq \sq{ \f{1}{4} + \f{
n_\ell \vep^2 } {2 \ln (\ze \de) } } \ee or \be \la{eqD2} \li |
\wh{p}_\ell - \f{1}{2} \ri | - \f{2\vep}{3} \leq - \sq{ \f{1}{4} +
\f{  n_\ell \vep^2 } {2 \ln (\ze \de) } }. \ee Since (\ref{eqD1})
implies either {\small $\wh{p}_\ell - \f{1}{2} + \f{2\vep}{3}  \geq
\f{4\vep}{3} + \sq{ \f{1}{4} + \f{  n_\ell \vep^2 } {2 \ln (\ze \de)
} } > \sq{ \f{1}{4} + \f{  n_\ell \vep^2 } {2 \ln (\ze \de) } }$} or
{\small $\wh{p}_\ell - \f{1}{2} + \f{2\vep}{3} \leq - \sq{ \f{1}{4}
+ \f{  n_\ell \vep^2 } {2 \ln (\ze \de) } }$}, it must be true that
(\ref{eqD1}) implies (\ref{eqDD}).  On the other hand, (\ref{eqD2})
also implies (\ref{eqDD}) because (\ref{eqD2}) implies {\small $\sq{
\f{1}{4} + \f{  n_\ell \vep^2 } {2 \ln (\ze \de) } } \leq
\wh{p}_\ell - \f{1}{2}  + \f{2\vep}{3}$}. Hence, we have established
(\ref{eqDD}) based on (\ref{eqD}).

 Since $- \f{1}{2} <  \wh{p}_\ell - \f{1}{2} + \f{2 \vep}{3} \leq p - \vep - \f{1}{2} + \f{2 \vep}{3} <
\f{1}{2}$,  we have $\f{1}{4} - \li ( \wh{p}_\ell - \f{1}{2} + \f{2
\vep}{3} \ri )^2 > 0$ and, by virtue of (\ref{eqDD}),
\[
\mcal{M} \li (\wh{p}_\ell, \wh{p}_\ell + \vep \ri ) = - \f{ \vep^2 }
{ 2 \li [ \f{1}{4} - \li ( \wh{p}_\ell - \f{1}{2} + \f{2 \vep}{3}
\ri )^2 \ri ] }  \leq \f{\ln (\ze \de)}{n_\ell}.
\]
Since $ \wh{p}_\ell \leq p - \vep$, we have  $0 \leq \wh{p}_\ell <
\wh{p}_\ell + \vep \leq p < 1$.  Hence, using  the fact that
$\mcal{M}(z, \mu)$ is monotonically decreasing with respect to $\mu
\in \li (z, 1 \ri )$ as asserted by Lemma \ref{Ddec}, we have
$\mcal{M} \li ( \wh{p}_\ell, p \ri ) \leq \mcal{M} \li (
\wh{p}_\ell, \wh{p}_\ell + \vep \ri ) \leq \f{\ln (\ze
\de)}{n_\ell}$.  The proof of the lemma is thus completed.  \epf

\beL \la{abs382}
{\small  $\li \{  \wh{\bs{p}}_\ell \geq
  p + \vep, \; \bs{D}_\ell = 1 \ri \}  \leu \li \{ \wh{\bs{p}}_\ell > p, \; \mcal{M} \li ( \wh{\bs{p}}_\ell, p \ri )
  \leq \f{\ln (\ze \de) }{n_\ell} \ri \}$} for $\ell = 1,
\cd, s$.

\eeL

\bpf

Since {\small $\{ \bs{D}_\ell = 1 \} = \li \{ \li ( \li |
\wh{\bs{p}}_\ell - \f{1}{2} \ri | - \f{2\vep}{3} \ri )^2 \geq
 \f{1}{4} + \f{  n_\ell \vep^2 } {2 \ln (\ze \de) } \ri \}$}, it suffices to show {\small
\[
\li \{ \wh{\bs{p}}_\ell \geq p + \vep,  \; \li ( \li |
\wh{\bs{p}}_\ell - \f{1}{2} \ri | - \f{2\vep}{3} \ri )^2 \geq
 \f{1}{4} + \f{  n_\ell \vep^2 } {2 \ln (\ze \de) } \ri \} \subseteq \li \{ \wh{\bs{p}}_\ell > p, \;
\mcal{M} \li ( \wh{\bs{p}}_\ell, p \ri ) \leq \f{\ln (\ze
\de)}{n_\ell} \ri \}. \]
 }
For this purpose, we let {\small $\om \in \li \{ \wh{\bs{p}}_\ell
\geq p + \vep, \; \li ( \li | \wh{\bs{p}}_\ell - \f{1}{2} \ri | -
\f{2\vep}{3} \ri )^2 \geq \f{1}{4} + \f{  n_\ell \vep^2 } {2 \ln
(\ze \de) } \ri \}, \; \wh{p}_\ell = \wh{\bs{p}}_\ell (\om)$} and
proceed to show $\wh{p}_\ell > p, \; \mcal{M} \li ( \wh{p}_\ell, p
\ri ) \leq \f{\ln (\ze \de)}{n_\ell}$. Clearly, $\wh{p}_\ell
> p$ follows immediately from $\wh{p}_\ell \geq p + \vep$.  To show
$\mcal{M} \li ( \wh{p}_\ell, p \ri ) \leq \f{\ln (\ze
\de)}{n_\ell}$, we need to establish \be \la{eqDD8} \li (
\wh{p}_\ell - \f{1}{2}  - \f{2\vep}{3} \ri )^2 \geq \f{1}{4} + \f{
n_\ell \vep^2 } {2 \ln (\ze \de) } \ee based on \be \la{eqD8}
 \li ( \li |
\wh{p}_\ell - \f{1}{2} \ri | - \f{2\vep}{3} \ri )^2 \geq
 \f{1}{4} + \f{  n_\ell \vep^2 } {2 \ln (\ze \de) }. \ee  It is obvious that (\ref{eqDD8}) holds if
{\small $\f{1}{4} + \f{  n_\ell \vep^2 } {2 \ln (\ze \de) } \leq
0$}.  It remains to show (\ref{eqDD8}) under the condition that
{\small $\f{1}{4} + \f{  n_\ell \vep^2 } {2 \ln (\ze \de) }
> 0$}.  Note that (\ref{eqD8}) implies either \be \la{eqD18} \li |
\wh{p}_\ell - \f{1}{2} \ri | - \f{2\vep}{3} \geq \sq{ \f{1}{4} + \f{
n_\ell \vep^2 } {2 \ln (\ze \de) } } \ee or \be \la{eqD28} \li |
\wh{p}_\ell - \f{1}{2} \ri | - \f{2\vep}{3} \leq - \sq{ \f{1}{4} +
\f{  n_\ell \vep^2 } {2 \ln (\ze \de) } }. \ee Since (\ref{eqD18})
implies either {\small $\wh{p}_\ell - \f{1}{2} - \f{2\vep}{3}  \leq
- \f{4\vep}{3} - \sq{ \f{1}{4} + \f{  n_\ell \vep^2 } {2 \ln (\ze
\de) } } < - \sq{ \f{1}{4} + \f{  n_\ell \vep^2 } {2 \ln (\ze \de) }
}$} or {\small $\wh{p}_\ell - \f{1}{2} - \f{2\vep}{3}  \geq  \sq{
\f{1}{4} + \f{  n_\ell \vep^2 } {2 \ln (\ze \de) } }$}, it must be
true that (\ref{eqD18}) implies (\ref{eqDD8}). On the other hand,
(\ref{eqD28}) also implies (\ref{eqDD8}) because (\ref{eqD28})
implies {\small $\wh{p}_\ell - \f{1}{2} - \f{2\vep}{3} \leq - \sq{
\f{1}{4} + \f{  n_\ell \vep^2 } {2 \ln (\ze \de) } }$}.  Hence, we
have established (\ref{eqDD8}) based on (\ref{eqD8}).

 Since $- \f{1}{2} < p + \vep - \f{1}{2} - \f{2 \vep}{3} \leq
\wh{p}_\ell - \f{1}{2} - \f{2 \vep}{3} \leq 1 - \f{1}{2} - \f{2
\vep}{3} < \f{1}{2}$,  we have $\f{1}{4} - \li ( \wh{p}_\ell -
\f{1}{2} - \f{2 \vep}{3} \ri )^2 > 0$ and, by virtue of
(\ref{eqDD8}),
\[
\mcal{M} \li (\wh{p}_\ell, \wh{p}_\ell - \vep \ri ) = - \f{ \vep^2 }
{ 2 \li [ \f{1}{4} - \li ( \wh{p}_\ell - \f{1}{2} - \f{2 \vep}{3}
\ri )^2 \ri ] }  \leq \f{\ln (\ze \de)}{n_\ell}.
\]
Since $ \wh{p}_\ell \geq p + \vep$, we have  $0 < p \leq \wh{p}_\ell
- \vep < \wh{p}_\ell \leq 1$.  Hence, using  the fact that
$\mcal{M}(z, \mu)$ is monotonically increasing with respect to $\mu
\in \li (0, z \ri )$ as asserted by Lemma \ref{Ddec}, we have
$\mcal{M} \li ( \wh{p}_\ell, p \ri ) \leq \mcal{M} \li (
\wh{p}_\ell, \wh{p}_\ell - \vep \ri ) \leq \f{\ln (\ze
\de)}{n_\ell}$.  The proof of the lemma is thus completed.

\epf

\bsk

Now we are in a position to prove Theorem 1.  Since $\f{9}{24 \vep -
16 \vep^2} > 1$ for any $\vep \in (0,\f{1}{2})$,  we have $\tau >
0$. Hence, the sequence of sample sizes $n_1, \cd, n_s$ is
well-defined.  By Lemma \ref{absD}, the sampling must stopped at
some stage with index $\ell \in \{1, \cd, s \}$. This shows that the
sampling scheme is well-defined. Noting that $\{ \mbf{n} = n_\ell \}
\subseteq \{ \bs{D}_\ell = 1 \}$ for $\ell = 1, \cd, s$, we have
\bel \Pr \{ | \wh{\bs{p}} - p | \geq \vep \} & = & \sum_{\ell = 1}^s
\{ \wh{\bs{p}}_\ell \leq p - \vep, \; \mbf{n} = n_\ell \} +
\sum_{\ell = 1}^s \{ \wh{\bs{p}}_\ell \geq p + \vep, \;
\mbf{n} = n_\ell \} \nonumber\\
& \leq & \sum_{\ell = 1}^s \{ \wh{\bs{p}}_\ell \leq p - \vep, \;
\bs{D}_\ell = 1 \} + \sum_{\ell = 1}^s \{ \wh{\bs{p}}_\ell \geq p +
\vep, \; \bs{D}_\ell = 1 \}. \la{com1abs} \eel By Lemmas
\ref{abs381} and \ref{lem2}, \be \la{com2abs}
 \sum_{\ell = 1}^s \{
\wh{\bs{p}}_\ell \leq p - \vep, \; \bs{D}_\ell = 1 \}  \leq
\sum_{\ell = 1}^s \li \{ \wh{\bs{p}}_\ell < p, \; \mcal{M} \li (
\wh{\bs{p}}_\ell, p \ri ) \leq  \f{\ln (\ze \de)}{n_\ell} \ri \}
\leq s \ze \de \leq (\tau + 1) \ze \de. \ee By Lemmas \ref{abs382}
and \ref{lem1}, \be \la{com3abs} \sum_{\ell = 1}^s \{
\wh{\bs{p}}_\ell \geq p + \vep, \; \bs{D}_\ell = 1 \} \leq
\sum_{\ell = 1}^s \li \{ \wh{\bs{p}}_\ell
> p, \; \mcal{M} \li ( \wh{\bs{p}}_\ell, p \ri ) \leq  \f{\ln (\ze \de)}{n_\ell} \ri \}
\leq s \ze \de \leq (\tau + 1) \ze \de. \ee Combining
(\ref{com1abs}), (\ref{com2abs}) and (\ref{com3abs}) yields $\Pr \{
| \wh{\bs{p}} - p | \geq \vep \} \leq 2 (\tau + 1) \ze \de$. Hence,
if we choose $\ze$ to be a positive number less than $\f{1}{2 (\tau
+ 1) }$, we have $\Pr \{  | \wh{\bs{p}} - p | < \vep \} > 1 - \de$.
This completes the proof of Theorem 1.

\sect{Proof of Theorem 2}

Throughout the proof of Theorem 2, we define
\[
\udl{\bs{p}}_\ell = \min \li \{ \wh{\bs{p}}_\ell - \vep_a, \;
\f{\wh{\bs{p}}_\ell}{1 + \vep_r} \ri \}, \qqu \ovl{\bs{p}}_\ell =
\max \li \{ \wh{\bs{p}}_\ell + \vep_a, \; \f{\wh{\bs{p}}_\ell}{1 -
\vep_r} \ri \}.
\]

By tedious computation, we can show the following lemma. \beL
\la{tran} For $\ell = 1, \cd, s$, {\small \be \la{trans1}
 \li \{ \wh{\bs{p}}_\ell
\geq \f{ 6(1 + \vep_r) (3 + \vep_r) \ln (\ze \de) } { 2 (3 +
\vep_r)^2 \ln (\ze \de) - 9 n_\ell \vep_r^2 } \ri \} = \li \{
\mcal{M} \li ( \wh{\bs{p}}_\ell, \f{\wh{\bs{p}}_\ell}{1 + \vep_r}
\ri ) \leq \f{ \ln (\ze \de) }{n_\ell} \ri \}, \ee \be \la{trans2}
 \li \{ \wh{\bs{p}}_\ell \geq \f{ 6(1 - \vep_r) (3 -
\vep_r) \ln (\ze \de) } { 2 (3 - \vep_r)^2 \ln (\ze \de) -  9 n_\ell
\vep_r^2} \ri \} = \li \{ \mcal{M} \li ( \wh{\bs{p}}_\ell,
\f{\wh{\bs{p}}_\ell}{1 - \vep_r} \ri ) \leq \f{ \ln (\ze \de)
}{n_\ell} \ri \}. \ee } \eeL

\beL

\la{bb3800}

{\small $\li \{ \udl{\bs{p}}_s \geq p \ri \} \subseteq \li \{
\wh{\bs{p}}_s > p, \; \mcal{M} \li ( \wh{\bs{p}}_s, p \ri ) \leq
\f{\ln (\ze \de)}{n_s} \ri \}$. }

\eeL

\bpf

To prove the lemma, we let $\om \in \{  \udl{\bs{p}}_s \geq p \}, \;
\wh{p}_s = \wh{\bs{p}}_s (\om), \; \udl{p}_s = \udl{\bs{p}}_s (\om)$
and proceed to show {\small $\wh{p}_s >  p, \; \mcal{M} \li (
\wh{p}_s, p \ri ) \leq \f{\ln (\ze \de)}{n_s}$}. Clearly, $\wh{p}_s
>  p$ follows immediately from $\udl{p}_s \geq p
> 0$.  To show {\small $\mcal{M} \li ( \wh{p}_s, p \ri ) \leq \f{\ln (\ze \de)}{n_s}$},
we shall first show {\small $\mcal{M} ( \wh{p}_s, \udl{p}_s ) \leq
\f{\ln (\ze \de)}{n_s}$}. For simplicity of notations, we denote
$p^\star = \f{\vep_a}{\vep_r}$. We need to consider three cases as
follows.

\bsk

Case (i):  $\wh{p}_s \leq p^\star - \vep_a$. In this case, {\small
\[ \mcal{M} ( \wh{p}_s, \udl{p}_s ) = \mcal{M} \li
(\wh{p}_s, \wh{p}_s - \vep_a \ri )
 <  \mcal{M} \li (\wh{p}_s, \wh{p}_s + \vep_a \ri ) \leq  \mcal{M} \li ( p^\star - \vep_a, p^\star \ri )  <
    \mcal{M} \li ( p^\star + \vep_a, p^\star \ri )
     \leq  \f{\ln (\ze \de)}{n_s}.
     \]}
Here the first inequality is due to $\vep_a < p + \vep_a \leq
\udl{p}_s + \vep_a = \wh{p}_s \leq p^\star - \vep_a < \f{1}{2}$ and
the fact that $\mcal{M}(z, z + \vep) > \mcal{M}(z, z - \vep)$ for
$\vep < z < \f{1}{2}$, which is asserted by Lemma \ref{lem888m}. The
second inequality is due to $\vep_a < p + \vep_a \leq \udl{p}_s +
\vep_a = \wh{p}_s < p^\star - \vep_a < \f{1}{2} - \vep_a$ and the
fact that $\mcal{M}(z, z + \vep)$ is monotonically increasing with
respect to $z \in (0, \f{1}{2} - \vep)$, which can be seen from
Lemma \ref{decrea}. The third inequality is due to $\vep_a < p^\star
< \f{1}{2}$ and the fact that $\mcal{M}(p + \vep, p)
> \mcal{M}(p - \vep, p)$ for $\vep < p < \f{1}{2}$ as a result of Lemma \ref{lemm22m}.  The last inequality
is due to the fact that {\small $n_s = \li \lc  \f{ \ln (\ze \de) }
{ \mcal{M} \li ( p^\star + \vep_a, p^\star \ri ) } \ri \rc \geq \f{
\ln (\ze \de) } { \mcal{M} \li ( p^\star + \vep_a, p^\star \ri ) }
$}, which follows from the definition of $n_s$.

\bsk

Case (ii):  $p^\star - \vep_a < \wh{p}_s < p^\star + \vep_a$. In
this case, {\small \[ \mcal{M} (\wh{p}_s,  \udl{p}_s )  = \mcal{M}
\li (\wh{p}_s, \wh{p}_s - \vep_a \ri )
 < \mcal{M} \li ( p^\star + \vep_a, p^\star + \vep_a - \vep_a
\ri )
 =  \mcal{M} \li ( p^\star + \vep_a, p^\star \ri)  \leq \f{\ln (\ze \de)}{n_s} \]} where the first inequality  is due to $
\vep_a < p + \vep_a \leq \udl{p}_s + \vep_a = \wh{p}_s < p^\star +
\vep_a < \f{1}{2} - \f{\vep_a}{3} + \vep_a$ and the fact that
$\mcal{M}(z, z - \vep)$ is monotonically increasing with respect to
$z \in (\vep, \f{1}{2} + \f{2 \vep}{3})$, which can be seen from
Lemma \ref{decrea}.

\bsk

Case (iii):  $\wh{p}_s \geq p^\star + \vep_a$.  In this case,
{\small
\[ \mcal{M} (\wh{p}_s, \udl{p}_s )
 =  \mcal{M}\li (\wh{p}_s, \f{\wh{p}_s}{1 + \vep_r} \ri )
 \leq  \mcal{M} \li (p^\star + \vep_a, \f{p^\star + \vep_a}{1 + \vep_r} \ri )
 =  \mcal{M} \li ( p^\star + \vep_a, p^\star \ri ) \leq  \f{\ln (\ze \de)}{n_s}.
 \]}
 where the first
inequality  is due to the fact that $\mcal{M}(z, z \sh (1 + \vep) )$
is monotonically decreasing with respect to $z \in (0, 1)$, which
can be seen from Lemma \ref{decrev}.

Therefore, we have shown {\small $\mcal{M} ( \wh{p}_s, \udl{p}_s  )
\leq \f{\ln (\ze \de)}{n_s}$} for all cases.   As a result of Lemma
\ref{Ddec}, $\mcal{M}(z, \mu)$ is monotonically increasing with
respect to $\mu \in (0, z)$.  By virtue of such monotonicity and the
fact that $\wh{p}_s \geq \udl{p}_s \geq p > 0$, we have {\small
$\mcal{M} \li ( \wh{p}_s, p \ri ) \leq \mcal{M} ( \wh{p}_s,
\udl{p}_s ) \leq \f{\ln (\ze \de)}{n_s}$}. This completes the proof
of the lemma.

 \epf

\beL

\la{bb380}

 {\small $\li \{  \ovl{\bs{p}}_s \leq p \ri \} \subseteq
\li \{ \wh{\bs{p}}_s <   p,  \; \mcal{M} \li ( \wh{\bs{p}}_s, p \ri
) \leq \f{\ln (\ze \de)}{n_s} \ri \}$. }

\eeL

\bpf

To prove the lemma, we let $\om \in \li \{  \ovl{\bs{p}}_s \leq p
\ri \}, \; \wh{p}_s = \wh{\bs{p}}_s (\om), \; \ovl{p}_s =
\ovl{\bs{p}}_s (\om)$ and proceed to show {\small $\wh{p}_s <   p,
\; \mcal{M} \li ( \wh{p}_s, p \ri ) \leq \f{\ln (\ze \de)}{n_s}$}.
Clearly, $\wh{p}_s <  p$ follows immediately from $\ovl{p}_s \leq p
< 1$.  To show {\small $\mcal{M} \li ( \wh{p}_s, p \ri ) \leq \f{\ln
(\ze \de)}{n_s}$},  we shall first show {\small $\mcal{M} \li (
\wh{p}_s, \ovl{p}_s \ri ) \leq \f{\ln (\ze \de)}{n_s}$} by
considering three cases as follows.

 \bsk

Case (i):  $\wh{p}_s \leq p^\star - \vep_a$. In this case, {\small
\[ \mcal{M} \li (\wh{p}_s,  \ovl{p}_s \ri )  =  \mcal{M} \li (
\wh{p}_s, \wh{p}_s + \vep_a \ri )
 \leq  \mcal{M} \li ( p^\star - \vep_a, p^\star - \vep_a + \vep_a \ri )
  =   \mcal{M} \li ( p^\star - \vep_a, p^\star \ri ) \nonumber\\
 <  \mcal{M} \li ( p^\star + \vep_a, p^\star \ri ) \leq  \f{\ln (\ze \de)}{n_s}.
 \]}
Here the first inequality is due to $0 \leq \wh{p}_s \leq p^\star -
\vep_a < \f{1}{2} - \vep_a$ and the fact that $\mcal{M}(z, z +
\vep)$ is monotonically increasing with respect to $z \in (0,
\f{1}{2} - \vep)$, which is asserted by Lemma \ref{decrea}.  The
second inequality is due to $\vep_a < p^\star < \f{1}{2}$ and the
fact that $\mcal{M}(p + \vep, p)
> \mcal{M}(p - \vep, p)$ for $\vep < p < \f{1}{2}$, which can be seen from Lemma \ref{lemm22m}.

\bsk

Case (ii):  $p^\star - \vep_a < \wh{p}_s < p^\star + \vep_a$. In
this case, {\small \[ \mcal{M} \li (\wh{p}_s,  \ovl{p}_s \ri ) =
\mcal{M} \li (\wh{p}_s, \f{\wh{p}_s}{1 - \vep_r} \ri ) <
 \mcal{M} \li ( p^\star - \vep_a, \f{p^\star - \vep_a}{1 - \vep_r} \ri )
  =  \mcal{M} \li (p^\star - \vep_a, p^\star \ri ) <
   \mcal{M} \li ( p^\star + \vep_a, p^\star \ri )
    \leq \f{\ln (\ze \de)}{n_s} \]}
 where the first inequality is due to $0 < p^\star - \vep_a < \wh{p}_s = (1 - \vep_r)
\ovl{p}_s \leq (1 - \vep_r) p < 1 - \vep_r$ and the fact that
$\mcal{M}(z, z \sh (1 - \vep) )$ is monotonically decreasing with
respect to $z \in (0, 1 - \vep)$, which is asserted by Lemma
\ref{decrev}.

\bsk

Case (iii):  $\wh{p}_s \geq p^\star + \vep_a$.  In this case,
{\small
\[ \mcal{M} \li (\wh{p}_s,  \ovl{p}_s \ri ) = \mcal{M} \li
(\wh{p}_s, \f{\wh{p}_s}{1 - \vep_r} \ri )
 <  \mcal{M} \li (\wh{p}_s, \f{\wh{p}_s}{1 + \vep_r} \ri )
 \leq  \mcal{M} \li (p^\star + \vep_a, \f{p^\star + \vep_a}{1 + \vep_r} \ri )
 =  \mcal{M} \li (p^\star + \vep_a, p^\star \ri ) \leq \f{\ln (\ze \de)}{n_s}.
 \]}
 Here the first inequality is due to $0 < \wh{p}_s = (1 - \vep_r) \ovl{p}_s \leq (1 -
\vep_r) p < 1 - \vep_r$ and the fact that $\mcal{M}(z, z \sh (1 +
\vep) ) > \mcal{M}(z, z \sh (1  - \vep) )$ for $0 < z < 1 - \vep$,
which can be seen from Lemma \ref{lem888m2B}. The second inequality
is due to $p^\star + \vep_a \leq \wh{p}_s$ and the fact that
$\mcal{M}(z, z \sh (1 + \vep) )$ is monotonically decreasing with
respect to $z \in (0, 1)$, which is asserted by Lemma \ref{decrev}.

Therefore, we have shown {\small $\mcal{M} \li ( \wh{p}_s, \ovl{p}_s
\ri ) \leq \f{\ln (\ze \de)}{n_s}$} for all cases.  As a result of
Lemma \ref{Ddec}, $\mcal{M}(z, \mu)$ is monotonically decreasing
with respect to $\mu \in (z, 1)$.  By virtue of such monotonicity
and the fact that $\wh{p}_s \leq \ovl{p}_s \leq p < 1$, we have
{\small $\mcal{M} \li ( \wh{p}_s, p \ri ) \leq \mcal{M} \li (
\wh{p}_s, \ovl{p}_s \ri ) \leq \f{\ln (\ze \de)}{n_s}$}. This
completes the proof of the lemma.

\epf

\beL

\la{bb381}

{\small $\li \{  \ovl{\bs{p}}_\ell \leq p, \; \bs{D}_\ell = 1 \ri \}
\subseteq \li \{ \wh{\bs{p}}_\ell <   p,  \; \mcal{M} \li (
\wh{\bs{p}}_\ell, p \ri ) \leq  \f{\ln (\ze \de)}{n_\ell} \ri \}$}
for $1 \leq \ell < s$.

\eeL

\bpf

To show the lemma, we let $\om \in \{  \ovl{\bs{p}}_\ell \leq p, \;
\bs{D}_\ell = 1 \}, \; \wh{p}_\ell = \wh{\bs{p}}_\ell (\om), \;
\ovl{p}_\ell = \ovl{\bs{p}}_\ell (\om)$ and proceed to show {\small
$\wh{p}_\ell <   p,  \; \mcal{M} \li ( \wh{p}_\ell, p \ri ) \leq
\f{\ln (\ze \de)}{n_\ell}$}.  Clearly, $\wh{p}_\ell < p$ follows
immediately from $\ovl{p}_\ell \leq p < 1$. To show {\small
$\mcal{M} \li ( \wh{p}_\ell, p \ri ) \leq \f{\ln (\ze
\de)}{n_\ell}$}, we shall first show {\small $\mcal{M} \li (
\wh{p}_\ell, \ovl{p}_\ell \ri ) \leq \f{\ln (\ze \de)}{n_\ell}$} by
considering three cases as follows.

Case (i): $\wh{p}_\ell  \leq \f{ \vep_a  } { \vep_r } - \vep_a$. In
this case, by the definition of the stopping rule, we have {\small
$\wh{p}_\ell \leq \f{1}{2} - \f{2\vep_a}{3} - \sq{ \f{1}{4} +  \f{
n_\ell \; \vep_a^2 } {2 \ln (\ze \de) } }$}, which implies {\small
$\f{1}{4} - \li ( \wh{p}_\ell - \f{1}{2} +  \f{2 \vep_a}{3} \ri )^2
\leq  - \f{  n_\ell \; \vep_a^2 } {2 \ln (\ze \de) }$}. Observing
that $\ovl{p}_\ell = \wh{p}_\ell + \vep_a \leq p < 1$, we have
\[
- \f{1}{2} < \wh{p}_\ell - \f{1}{2} +  \f{2 \vep_a}{3} =
\ovl{p}_\ell - \vep_a - \f{1}{2} +  \f{2 \vep_a}{3} \leq p - \vep_a
- \f{1}{2} +  \f{2 \vep_a}{3} < \f{1}{2}.
\]
Hence, $\f{1}{4} - \li ( \wh{p}_\ell - \f{1}{2} +  \f{2 \vep_a}{3}
\ri )^2 > 0$ and $\mcal{M} \li ( \wh{p}_\ell, \ovl{p}_\ell \ri ) = -
\f{\vep_a^2}{2 \li [ \f{1}{4} - \li ( \wh{p}_\ell - \f{1}{2} + \f{2
\vep_a}{3} \ri )^2 \ri ]} \leq \f{\ln (\ze \de)}{n_\ell}$.

Case (ii):  $\li | \wh{p}_\ell - \f{\vep_a}{\vep_r} \ri | < \vep_a$.
In this case, by the definition of the stopping rule and
(\ref{trans2}) of Lemma \ref{tran}, we have $\mcal{M} \li (
\wh{p}_\ell, \ovl{p}_\ell \ri ) \leq \f{\ln (\ze \de)}{n_\ell}$ with
$\ovl{p}_\ell = \f{\wh{p}_\ell}{1 - \vep_r}$.

Case (iii):  $\wh{p}_\ell  \geq \f{ \vep_a } { \vep_r } + \vep_a$.
In this case, we have $\ovl{p}_\ell = \f{\wh{p}_\ell}{1 - \vep_r}$
and {\small $\mcal{M} \li ( \wh{p}_\ell, \ovl{p}_\ell \ri ) =
\mcal{M} \li ( \wh{p}_\ell, \f{\wh{p}_\ell} {1 - \vep_r } \ri ) <
\mcal{M} \li ( \wh{p}_\ell, \f{\wh{p}_\ell} {1 + \vep_r } \ri ) \leq
\f{\ln (\ze \de)}{n_\ell}$}.   Here, the first inequality follows
from Lemma \ref{lem888m2B} and the fact that $\wh{p}_\ell = (1 -
\vep_r) \ovl{p}_\ell \leq (1 - \vep_r) p < 1 - \vep_r$.  The second
inequality follows from the definitions of the stopping rule and
(\ref{trans1}) of Lemma \ref{tran}.

Therefore, we have shown {\small $\mcal{M} \li ( \wh{p}_\ell,
\ovl{p}_\ell \ri ) \leq \f{\ln (\ze \de)}{n_\ell}$} for all three
cases.  As a result of Lemma \ref{Ddec}, $\mcal{M}(z, \mu)$ is
monotonically decreasing with respect to $\mu \in (z, 1)$.  By
virtue of such monotonicity and the fact that $0 < \wh{p}_\ell <
\ovl{p}_\ell \leq p < 1$, we have $\mcal{M} \li ( \wh{p}_\ell, p \ri
) \leq \mcal{M} \li ( \wh{p}_\ell,  \ovl{p}_\ell \ri ) \leq \f{\ln
(\ze \de)}{n_\ell}$.  This completes the proof of the lemma.

\epf

\beL

\la{bb382}

{\small $\li \{  \udl{\bs{p}}_\ell  \geq p, \; \bs{D}_\ell = 1 \ri
\} \subseteq \li \{ \wh{\bs{p}}_\ell > p, \; \mcal{M} \li (
\wh{\bs{p}}_\ell, p \ri ) \leq  \f{\ln (\ze \de)}{n_\ell} \ri \}$}
for $1 \leq \ell < s$.

\eeL

\bpf

To show the lemma, we let $\om \in  \{  \udl{\bs{p}}_\ell \geq p, \;
\bs{D}_\ell = 1 \}, \; \wh{p}_\ell = \wh{\bs{p}}_\ell (\om), \;
\udl{p}_\ell = \udl{\bs{p}}_\ell (\om)$ and proceed to show {\small
$\wh{p}_\ell >   p,  \; \mcal{M} \li ( \wh{p}_\ell, p \ri ) \leq
\f{\ln (\ze \de)}{n_\ell}$}.  Clearly, $\wh{p}_\ell > p$ follows
immediately from $\udl{p}_\ell \geq p > 0$. To show {\small
$\mcal{M} \li ( \wh{p}_\ell, p \ri ) \leq \f{\ln (\ze
\de)}{n_\ell}$}, we shall first show {\small $\mcal{M} (
\wh{p}_\ell, \udl{p}_\ell  ) \leq \f{\ln (\ze \de)}{n_\ell}$} by
considering three cases as follows.

Case (i):  $\wh{p}_\ell  \leq \f{ \vep_a  } { \vep_r } - \vep_a$. In
this case, we have $\udl{p}_\ell = \wh{p}_\ell - \vep_a$ and
\[
\vep_a < p + \vep_a \leq \udl{p}_\ell + \vep_a = \wh{p}_\ell \leq
\f{ \vep_a  } { \vep_r } - \vep_a \leq \f{1}{2} - \f{ 4 \vep_a } {
3},
\]
where the last inequality follows from the assumption about $\vep_a$
and $\vep_r$.   By virtue of the fact that $\vep_a < \wh{p}_\ell  <
\f{1}{2} < 1 - \vep_a$ and Lemma \ref{lem888m}, we have {\small
$\mcal{M} ( \wh{p}_\ell, \udl{p}_\ell
 ) = \mcal{M} \li ( \wh{p}_\ell, \wh{p}_\ell - \vep_a \ri ) <
\mcal{M} \li ( \wh{p}_\ell, \wh{p}_\ell + \vep_a \ri )$}.  Since $-
\f{1}{2} < \vep_a - \f{1}{2} + \f{ 2 \vep_a } { 3} < \wh{p}_\ell -
\f{1}{2} + \f{ 2 \vep_a } { 3} < \f{1}{2}$, we have $\f{1}{4} - \li
( \wh{p}_\ell - \f{1}{2} +  \f{2 \vep_a}{3} \ri )^2 > 0$. By the
definition of the stopping rule, we have {\small $\wh{p}_\ell \leq
\f{1}{2} - \f{2\vep_a}{3} - \sq{ \f{1}{4} +  \f{ n_\ell \; \vep_a^2
} {2 \ln (\ze \de) } }$}, which implies {\small $\f{1}{4} - \li (
\wh{p}_\ell - \f{1}{2} +  \f{2 \vep_a}{3} \ri )^2 \leq  - \f{ n_\ell
\; \vep_a^2 } {2 \ln (\ze \de) }$} and thus $\mcal{M} ( \wh{p}_\ell,
\udl{p}_\ell ) < \mcal{M} \li ( \wh{p}_\ell, \wh{p}_\ell + \vep_a
\ri ) \leq \f{\ln (\ze \de)}{n_\ell}$.

Case (ii):  $\li | \wh{p}_\ell - \f{\vep_a}{\vep_r} \ri | < \vep_a$.
In this case,  since $\udl{p}_\ell = \wh{p}_\ell - \vep_a \geq p >
0$, we have \[ - \f{1}{2} < p + \vep_a - \f{1}{2} - \f{2 \vep_a}{3}
\leq \udl{p}_\ell + \vep_a - \f{1}{2} -  \f{2 \vep_a}{3} =
\wh{p}_\ell - \f{1}{2} -  \f{2 \vep_a}{3} < \f{1}{2},
\]
which implies $\f{1}{4} - \li ( \wh{p}_\ell - \f{1}{2} - \f{2
\vep_a}{3} \ri )^2
> 0$. By the definition of the stopping rule, we have $\wh{p}_\ell
\leq \f{1}{2} + \f{2 \vep_a}{3} - \sq{ \f{1}{4} + \f{ n_\ell \;
\vep_a^2 } {2 \ln (\ze \de) } }$, which implies {\small $\f{1}{4} -
\li ( \wh{p}_\ell - \f{1}{2} -  \f{2 \vep_a}{3} \ri )^2 \leq -  \f{
n_\ell \; \vep_a^2 } {2 \ln (\ze \de) }$}.  It follows that {\small
$\mcal{M} ( \wh{p}_\ell, \udl{p}_\ell ) = \mcal{M} \li (
\wh{p}_\ell, \wh{p}_\ell - \vep_a \ri ) \leq \f{\ln (\ze
\de)}{n_\ell}$} because $\f{1}{4} - \li ( \wh{p}_\ell - \f{1}{2} -
\f{2 \vep_a}{3} \ri )^2 > 0$.

Case (iii):  $\wh{p}_\ell  \geq \f{ \vep_a } { \vep_r } + \vep_a$.
In this case, by the definition of the stopping rule and
(\ref{trans1}) of Lemma \ref{tran}, we have {\small $\mcal{M} (
\wh{p}_\ell, \udl{p}_\ell ) \leq \f{\ln (\ze \de)}{n_\ell}$} with
$\udl{p}_\ell = \f{\wh{p}_\ell} {1 + \vep_r }$.

Therefore, we have shown {\small $\mcal{M} ( \wh{p}_\ell,
\udl{p}_\ell ) \leq \f{\ln (\ze \de)}{n_\ell}$} for all three cases.
As a result of Lemma \ref{Ddec}, $\mcal{M}(z, \mu)$ is monotonically
increasing with respect to $\mu \in (0, z)$.  By virtue of such
monotonicity and the fact that $0 < p \leq \udl{p}_\ell <
\wh{p}_\ell < 1$,  we have $\mcal{M} \li ( \wh{p}_\ell, p \ri ) \leq
\mcal{M} ( \wh{p}_\ell, \udl{p}_\ell ) \leq \f{\ln (\ze
\de)}{n_\ell}$.  This completes the proof of the lemma.  \epf

\bsk

Now we are in a position to prove Theorem 2.  As a direct
consequence of the assumption that $0 < \vep_a < \f{3}{8}$ and
{\small $\f{6 \vep_a}{3 - 2 \vep_a } < \vep_r < 1$}, we have {\small
$\f{3}{2} \li ( \f{1}{\vep_a} - \f{1}{\vep_r} - \f{1}{3} \ri )
> 1$}, which implies that $\tau > 0$.  This shows that the sequence
of sample sizes $n_1, \cd, n_s$ is well-defined and it follows that
the sampling scheme is well-defined.  Invoking the definitions of
$\udl{\bs{p}}_\ell, \; \ovl{\bs{p}}_\ell$ and noting that $\{
\mbf{n} = n_\ell \} \subseteq \{ \bs{D}_\ell = 1 \}$ for $\ell = 1,
\cd, s$, we have \bel \Pr \{  | \wh{\bs{p}} - p | \geq \vep_a, \;  |
\wh{\bs{p}} - p | \geq \vep_r p  \} & = & \sum_{\ell = 1}^s \Pr \{ |
\wh{\bs{p}}_\ell - p | \geq \vep_a, \;  | \wh{\bs{p}}_\ell - p |
\geq \vep_r p, \;
\mbf{n} = n_\ell \} \nonumber\\
& = & \sum_{\ell = 1}^s \{ \ovl{\bs{p}}_\ell \leq p, \; \mbf{n} =
n_\ell \} + \sum_{\ell = 1}^s \{ \udl{\bs{p}}_\ell \geq p, \;
\mbf{n} = n_\ell \} \nonumber\\
& \leq & \sum_{\ell = 1}^s \{ \ovl{\bs{p}}_\ell \leq p, \;
\bs{D}_\ell = 1 \} + \sum_{\ell = 1}^s \{ \udl{\bs{p}}_\ell \geq p,
\; \bs{D}_\ell = 1 \}. \la{com1} \eel By Lemmas \ref{bb381},
\ref{bb380} and \ref{lem2}, \be \la{com2}
 \sum_{\ell = 1}^s \{
\ovl{\bs{p}}_\ell \leq p, \; \bs{D}_\ell = 1 \}  \leq  \sum_{\ell =
1}^s \li \{ \wh{\bs{p}}_\ell < p, \; \mcal{M} \li (
\wh{\bs{p}}_\ell, p \ri ) \leq  \f{\ln (\ze \de)}{n_\ell} \ri \}
\leq s \ze \de \leq (\tau + 1) \ze \de. \ee By Lemmas \ref{bb382},
\ref{bb3800} and \ref{lem1}, \be \la{com3} \sum_{\ell = 1}^s \{
\udl{\bs{p}}_\ell \geq p, \; \bs{D}_\ell = 1 \} \leq \sum_{\ell =
1}^s \li \{ \wh{\bs{p}}_\ell
> p, \; \mcal{M} \li ( \wh{\bs{p}}_\ell, p \ri ) \leq  \f{\ln (\ze \de)}{n_\ell} \ri \}
\leq s \ze \de \leq (\tau + 1) \ze \de. \ee Combining (\ref{com1}),
(\ref{com2}) and (\ref{com3}) yields $\Pr \{  | \wh{\bs{p}} - p |
\geq \vep_a, \;  | \wh{\bs{p}} - p | \geq \vep_r p  \} \leq 2 (\tau
+ 1) \ze \de$. Hence, if we choose $\ze$ to be a positive number
less than $\f{1}{2 (\tau + 1) }$, we have $\Pr \{  | \wh{\bs{p}} - p
| < \vep_a \; \tx{or} \; | \wh{\bs{p}} - p | < \vep_r p  \} > 1 -
\de$. This completes the proof of Theorem 2.

\sect{Proof of Theorem 3}

In the course of proving Theorem 3, we need to use the following
lemma regarding inverse binomial sampling, which has been
established by Chen in \cite{Chen_EST}.

\beL \la{gen_H} Let $X_1, X_2, \cd$ be a sequence of i.i.d.
Bernoulli random variables such that $\Pr \{ X_i = 1 \} = 1 - \Pr \{
X_i = 0 \} = p \in (0, 1)$ for $i = 1, 2, \cd$. Let $\bs{n}$ be the
minimum integer such that $\sum_{i=1}^{\bs{n}} X_i = \ga$ where
$\ga$ is a positive integer.  Then, for any $\al > 0$, {\small
\[
 \Pr \li \{
\f{\ga}{\bs{n}} \leq p, \; \mscr{M}_{\mrm{I}} \li ( \f{\ga}{\bs{n}},
p \ri ) \leq \f{ \ln \al } { \ga } \ri \} \leq \al, \qqu
 \Pr \li \{ \f{\ga}{\bs{n}} \geq p, \;
 \mscr{M}_{\mrm{I}} \li ( \f{\ga}{\bs{n}}, p \ri ) \leq \f{ \ln \al } {
\ga } \ri \} \leq \al \]} where {\small \[ \mscr{M}_{\mrm{I}}
(z,\mu) = \bec  \ln \f{\mu}{z} + \li ( \f{1}{z} - 1 \ri ) \ln \f{1 -
\mu}{1 - z} &
\tx{for} \; z \in (0,1) \; \tx{and} \; \mu \in (0, 1),\\
\ln \mu &  \tx{for} \; z = 1 \; \tx{and} \; \mu \in (0, 1),\\
- \iy & \tx{for} \; z = 0 \; \tx{and} \; \mu \in (0, 1). \eec
\]}
\eeL

\beL \la{gen} Let $X_1, X_2, \cd$ be a sequence of i.i.d. Bernoulli
random variables such that $\Pr \{ X_i = 1 \} = 1 - \Pr \{ X_i = 0
\} = p \in (0, 1)$ for $i = 1, 2, \cd$. Let $\bs{n}$ be the minimum
integer such that $\sum_{i=1}^{\bs{n}} X_i = \ga$ where $\ga$ is a
positive integer.  Then, for any $\al > 0$, {\small \be \la{bb2} \Pr
\li \{ \f{\ga}{\bs{n}} \leq p, \; \mcal{M}_{\mrm{I}} \li (
\f{\ga}{\bs{n}}, p \ri ) \leq \f{ \ln \al } { \ga } \ri \} \leq \al,
\ee \be \la{bb1}
 \Pr \li \{ \f{\ga}{\bs{n}} \geq p, \;
 \mcal{M}_{\mrm{I}} \li ( \f{\ga}{\bs{n}}, p \ri ) \leq \f{ \ln \al } {
\ga } \ri \} \leq \al \ee} where $\mcal{M}_{\mrm{I}} (z, \mu) = \f{
\mcal{M}(z, \mu) }{z}$ for $0 < z \leq 1$ and $0 < \mu < 1$.

\eeL

\bpf

By Massart's inequality (i.e., Theorem 2 at page 1271 of
\cite{Massart:90}), we have $\mscr{M}_{\mrm{I}} (z, p) <
\mcal{M}_{\mrm{I}} (z, p)$ for any $z \in (0, p)$.  By virtue of
this fact and Lemma \ref{gen_H}, we have
\[
\Pr \li \{ \f{\ga}{\bs{n}} \leq p, \; \mcal{M}_{\mrm{I}} \li (
\f{\ga}{\bs{n}}, p \ri ) \leq \f{ \ln \al } { \ga } \ri \} \leq \Pr
\li \{ \f{\ga}{\bs{n}} \leq p, \; \mscr{M}_{\mrm{I}} \li (
\f{\ga}{\bs{n}}, p \ri ) \leq \f{ \ln \al } { \ga } \ri \} \leq \al,
\]
\[
\Pr \li \{ \f{\ga}{\bs{n}} \geq p, \; \mcal{M}_{\mrm{I}} \li (
\f{\ga}{\bs{n}}, p \ri ) \leq \f{ \ln \al } { \ga } \ri \} \leq \Pr
\li \{ \f{\ga}{\bs{n}} \geq p, \; \mscr{M}_{\mrm{I}} \li (
\f{\ga}{\bs{n}}, p \ri ) \leq \f{ \ln \al } { \ga } \ri \} \leq \al.
\]
This completes the proof of the lemma.

\epf

In the sequel, we define random variables $\bs{D}_\ell, \; \ell = 1,
\cd, s$ such that $\bs{D}_\ell = 1$ if {\small $\ga_\ell \geq \f{ 6
\mathbf{n}_\ell (1 + \vep) (3 + \vep) \ln (\ze \de) } { 2 (3 +
\vep)^2 \ln (\ze \de) -  9 \vep^2 \mathbf{n}_\ell}$} and
$\bs{D}_\ell = 0$ otherwise.  Then, the stopping rule can be
restated as ``sampling is continued until $\bs{D}_\ell = 1$ for some
$\ell \in \{1, \cd, s \}$''.  For simplicity of notations, we also
define $\bs{D}_0 = 0$.

By tedious computation, we can show the following lemma.

\beL \la{lem199}
 $\li \{  \bs{D}_\ell = 1 \ri \} = \li \{
\mcal{M}_{\mrm{I}} \li ( \wh{\bs{p}}_\ell, \f{\wh{\bs{p}}_\ell}{1 +
\vep} \ri ) \leq \f{\ln (\ze \de) }{\ga_\ell} \ri \}$ for $\ell = 1,
\cd, s$. \eeL

\beL \la{pass_a} $\li \{  \wh{\bs{p}}_\ell \leq p (1 - \vep), \;
\bs{D}_\ell = 1 \ri \} \subseteq \li \{  \wh{\bs{p}}_\ell < p, \;
\mcal{M}_{\mrm{I}} \li ( \wh{\bs{p}}_\ell,  p \ri ) \leq \f{\ln (\ze
\de) }{\ga_\ell} \ri \}$ for $\ell = 1, \cd, s$. \eeL

\bpf

Let {\small $\om \in \li \{  \wh{\bs{p}}_\ell \leq p (1 - \vep), \;
\bs{D}_\ell = 1 \ri \}$} and {\small $\wh{p}_\ell = \wh{\bs{p}}_\ell
(\om)$}.  To show the lemma, it suffices to show $\wh{p}_\ell < p$
and {\small $\mcal{M}_{\mrm{I}} \li ( \wh{p}_\ell, p \ri ) \leq
\f{\ln (\ze \de) }{\ga_\ell}$}.  By Lemma \ref{lem199}, {\small \[
\li \{  \wh{\bs{p}}_\ell \leq p (1 - \vep), \; \bs{D}_\ell = 1 \ri
\} = \li \{ \wh{\bs{p}}_\ell \leq p (1 - \vep), \;
\mcal{M}_{\mrm{I}} \li ( \wh{\bs{p}}_\ell , \f{\wh{\bs{p}}_\ell}{1 +
\vep} \ri ) \leq \f{\ln (\ze \de) } { \ga_\ell } \ri \}
\]}
which implies $\wh{p}_\ell \leq p (1 - \vep)$ and {\small
$\mcal{M}_{\mrm{I}} \li ( \wh{p}_\ell , \f{\wh{p}_\ell}{1 + \vep}
\ri ) \leq \f{\ln (\ze \de)} { \ga_\ell }$}. Clearly, $\wh{p}_\ell
\leq p (1 - \vep)$ implies $\wh{p}_\ell < p$. To show {\small
$\mcal{M}_{\mrm{I}} \li ( \wh{p}_\ell, p \ri ) \leq \f{\ln (\ze \de)
}{\ga_\ell}$}, we shall consider two cases as follows:

In the case $\wh{p}_\ell = 0$, we have {\small $\mcal{M}_{\mrm{I}}
\li ( \wh{p}_\ell, p \ri ) = - \iy < \f{\ln (\ze \de) }{\ga_\ell}$}.

In the case of $\wh{p}_\ell > 0$,  we have  $0 < \wh{p}_\ell \leq p
(1 - \vep) < 1 - \vep$. Since {\small \[ \mcal{M}_{\mrm{I}} \li ( z,
\f{z}{1 + \vep} \ri ) - \mcal{M}_{\mrm{I}} \li ( z, \f{z}{1 - \vep}
\ri ) = \f{ 2 \vep^3 (2 - z) } {3 \li ( 1 + \f{\vep}{3} \ri ) \li [
1 - z + \vep \li ( 1 - \f{ z}{3} \ri ) \ri ] \li ( 1 - \f{\vep}{3}
\ri ) \li [ 1 - z - \vep \li ( 1 - \f{ z}{3} \ri ) \ri ]} > 0
\]}
for $0 < z < 1 - \vep$,  we have {\small $\mcal{M}_{\mrm{I}} \li (
\wh{p}_\ell , \f{\wh{p}_\ell}{1 - \vep} \ri ) < \mcal{M}_{\mrm{I}}
\li ( \wh{p}_\ell , \f{\wh{p}_\ell}{1 + \vep} \ri ) \leq \f{\ln (\ze
\de)} { \ga_\ell }$}.  Note that {\small \[ \f{\pa
\mcal{M}_{\mrm{I}} (z,\mu)}{\pa \mu}  = \f{(z - \mu) \li [ \mu(1 -
z) + z ( 1 - \mu) + z ( 1 -z) \ri ] } { 3 z \li [ \li ( \f{2\mu}{3}
+ \f{z}{3} \ri ) \li ( 1 - \f{2\mu}{3} - \f{z}{3} \ri ) \ri ]^2 },
\]} from which it can be seen that $\mcal{M}_{\mrm{I}} (z, \mu)$ is
monotonically decreasing with respect to $\mu \in (z, 1)$.  By
virtue of such monotonicity and the fact that {\small $0 <
\wh{p}_\ell < \f{\wh{p}_\ell}{1 - \vep} \leq p < 1$}, we have
{\small $\mcal{M}_{\mrm{I}} \li ( \wh{p}_\ell, p \ri ) \leq
\mcal{M}_{\mrm{I}} \li ( \wh{p}_\ell , \f{\wh{p}_\ell}{1 - \vep} \ri
) < \f{\ln (\ze \de) }{\ga_\ell}$}.  This completes the proof of the
lemma.

\epf

\beL \la{pass_b} $\li \{  \wh{\bs{p}}_\ell \geq p (1 + \vep), \;
\bs{D}_\ell = 1 \ri \} \subseteq \li \{ \wh{\bs{p}}_\ell > p, \;
\mcal{M}_{\mrm{I}} \li ( \wh{\bs{p}}_\ell,  p \ri ) \leq  \f{\ln
(\ze \de) }{\ga_\ell} \ri \}$ for $\ell = 1, \cd, s$. \eeL

\bpf

Let {\small $\om \in \li \{  \wh{\bs{p}}_\ell \geq p (1 + \vep), \;
\bs{D}_\ell = 1 \ri \}$} and {\small $\wh{p}_\ell = \wh{\bs{p}}_\ell
(\om)$}.  To show the lemma, it suffices to show $\wh{p}_\ell > p$
and {\small $\mcal{M}_{\mrm{I}} \li ( \wh{p}_\ell, p \ri ) \leq
\f{\ln (\ze \de) }{\ga_\ell}$}.  By Lemma \ref{lem199}, {\small \[
\li \{  \wh{\bs{p}}_\ell \geq p (1 + \vep), \; \bs{D}_\ell = 1 \ri
\} = \li \{ \wh{\bs{p}}_\ell \geq p (1 + \vep), \;
\mcal{M}_{\mrm{I}} \li ( \wh{\bs{p}}_\ell , \f{\wh{\bs{p}}_\ell}{1 +
\vep} \ri ) \leq \f{\ln (\ze \de) } { \ga_\ell } \ri \}
\]}
which implies $\wh{p}_\ell \geq p (1 + \vep)$ and {\small
$\mcal{M}_{\mrm{I}} \li ( \wh{p}_\ell , \f{\wh{p}_\ell}{1 + \vep}
\ri ) \leq \f{\ln (\ze \de)} { \ga_\ell }$}. Clearly, $\wh{p}_\ell
\geq p (1 + \vep)$ implies $\wh{p}_\ell > p$.  Since $1 \geq
\wh{p}_\ell \geq p (1 + \vep)$, we have {\small $0 < p \leq \f{
\wh{p}_\ell } { 1 + \vep } < \wh{p}_\ell \leq 1$}.  Noting that
{\small $\f{ \pa \mcal{M}_{\mrm{I}} (z, \mu) } { \pa \mu } = \f{(z -
\mu) \li [ \mu(1 - z) + z ( 1 - \mu) + z ( 1 -z) \ri ] } { 3 z \li [
\li ( \f{2\mu}{3} + \f{z}{3} \ri ) \li ( 1 - \f{2\mu}{3} - \f{z}{3}
\ri ) \ri ]^2 } > 0$} for $ 0 < \mu < z < 1$, we have {\small
$\mcal{M}_{\mrm{I}} \li ( \wh{p}_\ell, p \ri ) \leq
\mcal{M}_{\mrm{I}} \li ( \wh{p}_\ell , \f{\wh{p}_\ell}{1 + \vep} \ri
) \leq \f{\ln (\ze \de) }{\ga_\ell}$}.  This completes the proof of
the lemma.

\epf

\beL \la{well}
 $\bs{D}_s = 1$.
\eeL

\bpf To show $\bs{D}_s = 1$, it suffices to show {\small
$\mcal{M}_{\mrm{I}} \li (z, \f{z}{1 + \vep} \ri ) \leq \f{ \ln (\ze
\de)  } { \ga_s }$} for any $z \in (0, 1]$.  This is because $0 <
\wh{\bs{p}}_s (\om) \leq 1$ for any $\om \in \Om$ and {\small $\{
\bs{D}_s = 1 \} = \li \{ \mcal{M}_{\mrm{I}} \li ( \wh{\bs{p}}_s ,
\f{\wh{\bs{p}}_s}{1 + \vep} \ri ) \leq \f{\ln (\ze \de) }{\ga_s} \ri
\}$} as asserted by Lemma \ref{lem199}.

By the definition of sample sizes, we have {\small $\ga_s = \li \lc
\f{ \ln (\ze \de)}{ -  \vep^2 \li [ 2 \li ( 1 + \f{\vep}{3} \ri ) (
1 + \vep ) \ri ]^{-1} } \ri \rc \geq \f{ \ln (\ze \de)}{ - \vep^2
\li [ 2 \li ( 1 + \f{\vep}{3} \ri ) ( 1 + \vep ) \ri ]^{-1} }$}.
Since {\small $\lim_{z \to 0} \mcal{M}_{\mrm{I}} \li (z, \f{z}{1 +
\vep} \ri )   = -  \vep^2 \li [ 2 \li ( 1 + \f{\vep}{3} \ri ) ( 1 +
\vep ) \ri ]^{-1}  < 0$}, we have {\small $\lim_{z \to 0}
\mcal{M}_{\mrm{I}} \li (z, \f{z}{1 + \vep} \ri ) \leq \f{ \ln (\ze
\de)  } { \ga_s }$}.

Note that {\small $\mcal{M}_{\mrm{I}} \li ( z, \f{z}{1 + \vep} \ri )
= - \f{ \vep^2}{2 \li ( 1 + \f{\vep}{3} \ri ) \li [ 1 + \vep - (1 -
\f{\vep}{3}) z \ri ] }$}, from which it can be seen that {\small
$\mcal{M}_{\mrm{I}} \li (z, \f{z}{1 + \vep} \ri )$} is monotonically
decreasing with respect to $z \in (0, 1)$. Hence, {\small
$\mcal{M}_{\mrm{I}} \li (z, \f{z}{1 + \vep} \ri ) < \lim_{z \to 0}
\mcal{M}_{\mrm{I}} \li (z, \f{z}{1 + \vep} \ri ) \leq \f{ \ln (\ze
\de)  } { \ga_s }$} for any $z \in (0, 1)$. Since {\small
$\mcal{M}_{\mrm{I}} \li (z, \f{z}{1 + \vep} \ri )$} is a continuous
function with respect to $z \in (0, 1)$ and {\small
$\mcal{M}_{\mrm{I}} \li (1, \f{1}{1 + \vep} \ri ) = \lim_{z \to 1}
\mcal{M}_{\mrm{I}} \li (z, \f{z}{1 + \vep} \ri )$}, it must be true
that {\small $\mcal{M}_{\mrm{I}} \li (1, \f{1}{1 + \vep} \ri ) \leq
\f{ \ln (\ze \de)  } { \ga_s }$}. This completes the proof of the
lemma.

\epf

\beL \la{revsuf1}

$\Pr  \{  \wh{\bs{p}} \leq p (1  - \vep )  \} \leq \sum_{\ell = 1}^s
\Pr \li \{ \wh{\bs{p}}_\ell \leq p (1 - \vep), \; \bs{D}_{\ell - 1}
= 0, \; \bs{D}_\ell = 1  \ri \}  \leq (\tau + 1) \ze \de$ for any $p
\in (0, 1)$.

 \eeL

\bpf

By Lemma \ref{well}, the sampling must stop at some stage with index
$\ell \in \{1, \cd, s\}$.  This implies that the stopping rule is
well-defined.  Let $\bs{\ga} = \sum_{i = 1}^{ \mathbf{n} } X_i$.
Then, we can write {\small $\Pr \{ \wh{\bs{p}} \leq p(1 - \vep) \} =
\sum_{\ell = 1}^s \Pr \{ \wh{\bs{p}}_\ell \leq p (1 - \vep ), \;
\bs{\ga} = \ga_\ell \}$}. By the definition of the stopping rule, we
have $ \{ \bs{\ga} = \ga_\ell \} \subseteq \{ \bs{D}_{\ell - 1}  =
0, \; \bs{D}_\ell = 1 \}$.  Hence, {\small \bel \Pr \li \{
\wh{\bs{p}} \leq  p(1 - \vep) \ri \} \leq \sum_{\ell = 1}^s  \Pr \li
\{ \wh{\bs{p}}_\ell \leq p (1 - \vep), \; \bs{D}_{\ell - 1}  = 0, \;
\bs{D}_\ell = 1 \ri \} \leq \sum_{\ell = 1}^s  \Pr \li \{
\wh{\bs{p}}_\ell \leq p (1 - \vep), \; \bs{D}_\ell = 1 \ri \}. \la
{summ1} \eel} Applying Lemma \ref{pass_a} and (\ref{bb2}) of Lemma
\ref{gen}, we have {\small \be \la{parta1}
 \sum_{\ell = 1}^s  \Pr \li \{
\wh{\bs{p}}_\ell \leq p (1 - \vep), \;  \bs{D}_\ell = 1 \ri \} \leq
\sum_{\ell = 1}^s \Pr \li \{ \wh{\bs{p}}_\ell < p, \;
\mcal{M}_{\mrm{I}} \li ( \wh{\bs{p}}_\ell, p \ri ) \leq \f{\ln (\ze
\de) }{\ga_\ell} \ri \} \leq s \ze \de \leq (\tau + 1) \ze \de. \ee}
Finally, the lemma can be established by combining (\ref{summ1}) and
(\ref{parta1}).

\epf

\beL \la{revsuf2}

$\Pr  \{  \wh{\bs{p}} \geq p (1  + \vep ) \} \leq \sum_{\ell = 1}^s
\Pr \li \{ \wh{\bs{p}}_\ell \geq p (1 + \vep), \; \bs{D}_{\ell - 1}
= 0, \; \bs{D}_\ell = 1 \ri \}  \leq (\tau + 1) \ze \de$ for any $p
\in (0, 1)$.

 \eeL

\bpf

Note that {\small \bel \Pr \li \{ \wh{\bs{p}} \geq  p(1 + \vep) \ri
\} \leq \sum_{\ell = 1}^s  \Pr \li \{ \wh{\bs{p}}_\ell \geq p (1 +
\vep), \; \bs{D}_{\ell - 1}  = 0, \; \bs{D}_\ell = 1 \ri \} \leq
\sum_{\ell = 1}^s  \Pr \li \{ \wh{\bs{p}}_\ell \geq p (1 + \vep), \;
\bs{D}_\ell = 1 \ri \}. \la {summ2} \eel} Applying Lemma
\ref{pass_b} and (\ref{bb1}) of Lemma \ref{gen}, we have {\small \be
\la{partb1}
 \sum_{\ell = 1}^s  \Pr \li \{
\wh{\bs{p}}_\ell \geq p (1 + \vep), \;  \bs{D}_\ell = 1 \ri \} \leq
\sum_{\ell = 1}^s \Pr \li \{ \wh{\bs{p}}_\ell > p, \;
\mcal{M}_{\mrm{I}} \li ( \wh{\bs{p}}_\ell, p \ri ) \leq \f{\ln (\ze
\de) }{\ga_\ell} \ri \} \leq s \ze \de \leq (\tau + 1) \ze \de. \ee}
Combining (\ref{summ2}) and (\ref{partb1}) proves the lemma.

\epf

\bsk

Finally, we are in a position to prove Theorem 3.  Noting that $\Pr
\{  | \wh{\bs{p}} - p | \geq \vep p  \} = \Pr  \{  \wh{\bs{p}} \leq
p (1  - \vep )  \} + \Pr  \{  \wh{\bs{p}} \geq p (1  + \vep ) \}$
and making use of Lemmas \ref{revsuf1} and \ref{revsuf2}, we have
$\Pr \{  | \wh{\bs{p}} - p | \geq \vep p  \} \leq (\tau + 1) \ze \de
+ (\tau + 1) \ze \de = 2 (\tau + 1) \ze \de$ for any $p \in (0, 1)$.
Hence, if we choose $\ze$ to be a positive number less than $\f{1}{2
(\tau + 1) }$, we have $\Pr \{  | \wh{\bs{p}} - p | \geq \vep p \} <
\de$ and thus $\Pr \{  | \wh{\bs{p}} - p | < \vep p \} > 1 - \de$
for any $p \in (0, 1)$.  This completes the proof of Theorem 3.

\end{document}